\documentclass[12pt]{article}
\usepackage{epsfig, epsf, graphicx}
\usepackage{pstricks, pst-node, psfrag}
\usepackage{amssymb,amsmath,bm}
\usepackage[colorlinks,citecolor=blue,urlcolor=blue]{hyperref}
\usepackage{color}
\usepackage{verbatim,enumerate}
\usepackage{rotating, lscape,natbib}
\usepackage{setspace}
\usepackage{subfigure}
\usepackage{tikz}
\usetikzlibrary{arrows,calc,tikzmark}
\usepackage{hyperref}
\usepackage{float}
\usepackage{caption}
\usepackage{natbib}
\usepackage{comment}
\usepackage{dsfont}
\usepackage{booktabs}
\usepackage{authblk} 
\usepackage{algorithm}
\usepackage{algpseudocode}

\makeatletter
\DeclareRobustCommand{\sqcdot}{\mathbin{\mathpalette\morphic@sqcdot\relax}}
\newcommand{\morphic@sqcdot}[2]{%
  \sbox\z@{$\m@th#1\centerdot$}%
  \ht\z@=.33333\ht\z@
  \vcenter{\box\z@}%
}
\makeatother

\numberwithin{equation}{section}

\newtheorem{Def}{Definition}

\newtheorem{rem}{Remark}

\begin{document}

\thispagestyle{empty} \baselineskip=28pt \vskip 5mm
\begin{center} {\Huge{\bf Recent Developments in Complex and Spatially Correlated Functional Data}}
	
\end{center}

\baselineskip=12pt \vskip 10mm

\begin{center}\large
Israel Mart\'inez-Hern\'andez$^{*}$\footnote[1]{\label{note1} 
\baselineskip=10pt Statistics Program,
King Abdullah University of Science and Technology,
Thuwal 23955-6900, Saudi Arabia.\\
\vspace{-.2cm}

 $^{*}$ Corresponding author. E-mail: israel.martinezhernandez@kaust.edu.sa \\
 
This publication is based on research supported by the King Abdullah University of Science and Technology (KAUST) Office of Sponsored Research (OSR) under Award No: OSR-2018-CRG7-3742.
 } and Marc G.~Genton$^{\ref{note1}}$
\end{center}

\baselineskip=17pt \vskip 10mm \centerline{\today} \vskip 15mm

\begin{center}
{\large{\bf Abstract}}
\end{center}
As high-dimensional and high-frequency data are being collected on a large scale, the development of new statistical models  is being pushed forward. Functional data analysis provides the required statistical methods to deal with large-scale and complex data by assuming that data are  continuous functions, e.g., realizations of a continuous process (curves) or continuous random field (surfaces), and that each curve or surface is considered as a single observation.  
Here, we provide an overview of functional data analysis when data are complex and spatially correlated. We provide definitions and estimators of the first and second moments of the corresponding functional random variable.
We present two main approaches: The first assumes that data are realizations of a functional random field, i.e., each observation is a curve with a  spatial component. We call them \textit{spatial functional data}. The second approach assumes that data are continuous deterministic fields observed over time. In this case, one observation is a surface or manifold, and we call them \textit{surface time series}. For these two approaches, we describe   
  software available for the statistical analysis. We also present a data illustration, using a high-resolution wind speed simulated dataset, as an example of the two approaches. 
The functional data approach offers a new paradigm of data analysis, where the continuous processes or random fields are considered as a single entity. We consider  this approach to be very valuable in the context of big data.

\baselineskip=14pt

\par\vfill\noindent
{\bf Keywords:} Functional data; Functional random field; Manifold data; Spatial functional data; Spatial statistics; Spatio-temporal statistics; Surface data. 
\par\medskip\noindent
{\bf Short title}: Complex and Spatially Correlated Functional Data

\clearpage\pagebreak\newpage \pagenumbering{arabic}
\baselineskip=26pt

\section{Introduction}\label{sec:intro}

The statistical analysis of large, complex, and high-dimensional data has become a significant challenging problem.  Due to the rapid development of complex, performant technologies, data can now be collected on a large scale, resulting in high-dimensional and high-frequency data, sometimes necessitating high-performance computing, which is often a limitation for practitioners; see \cite{Galeano2019} for a general view of data science and big data. Among various approaches, functional data analysis (FDA)  provides statistical methods to handle  large-scale and complex data \citep[][]{Chen2017,GIRALDO2018126}. For a general  introduction to FDA, the reader is referred to  \cite{Ramsay-Silverman2005}, \cite{FerratyVieu2006}, \cite{HorvathKokoszka2012}, and \cite{KokoszkaEtAl2017}. FDA assumes that observations (called functional data) have characteristics that vary along a continuum, e.g., curves or surfaces. Thus, FDA deals with data that are defined on a space that is intrinsically infinite-dimensional.

The approach of FDA has several advantages as a realization of a continuous process or a continuous random field can be  considered as functional data. In this case, the stationarity of the process (random field) is not needed since FDA treats the whole curve (surface) as a single entity. Thus, FDA is part of object data analysis \citep[][]{MENAFOGLIO2017401}.  FDA is useful when the number of variables, $p$, is bigger than the sample size, $n$: $p\gg n$. In particular, FDA can analyze longitudinal data. Smoothness is an important property of functional data, in contrast with multivariate data analysis, where smoothness has no meaning. Thus, FDA extracts additional information contained in a  continuous function or in its derivative. Although in practice, each continuous functions, say $y_{i}(v)$,  is observed on a finite set of points, the continuity is obtained with smoothing techniques. In the process of smoothing, FDA does not require observed $y_{i}(v)$ on a regular grid, that is, a sample of $y_{i}(v)$ and $y_{j}( v)$ can be collected on a different set of points; $\mathbf{v}_{i}= \{v_{1_{i}}, \ldots, v_{n_{i}} \}$    and $\mathbf{v}_{j} = \{v_{1_{j}}, \ldots, v_{n_{j}} \}$, respectively. In general, the methods of FDA are essentially nonparametric and can model complex and spatially correlated data.  

Functional data can also have a spatial component, because data are collected somewhere at some time \cite[][p.15]{haining2003}. If the functional data are curves with a spatial component, we call them \textit{spatial functional data} (SFD). Thus, a dataset of SFD has the form $y(\mathbf{s}_{1}; v), \ldots, y(\mathbf{s}_{n}; v)$, where $\mathbf{s}_{i}\in D$ represents the locations in a given region $D$, and $v$ represents the continuous parameter of the functional data. For example, $y(\mathbf{s}_{i}; v)$ can be the daily wind speed observed at location  $\mathbf{s}_{i}$, with $v$ indicating the time within a day, or $y(\mathbf{s}_{i}; v)$ can be the spectrum of brain activity at location $\mathbf{s}_{i}$, with $v$ representing the frequency. If the continuous parameter $v$ represents time, as for the example of wind speed, then SFD can be related to spatio-temporal data where the temporal dependence is captured through the continuity of the curve. 

The combination of FDA and spatial statistics provides a  powerful tool to deal with complex and large spatial data. This combination is attracting   interest, and much research is focused on this topic. \cite{NERINI2010} proposed a spatial functional linear model, and they analyzed data in Oceanography. \cite{ZhouEtAl2010} proposed mixed effects models for spatially correlated hierarchical functional data. \cite{RUIZMEDINA2011292} extended the spatial autoregressive processes and the spatial moving average processes to the Hilbert space. \cite{GiraldoEtAl2012} proposed a methodology for clustering spatially correlated functional data; see also, \cite{HuijingEtAl2012} and \cite{Romano2017}.  \cite{StaicuEtAl2010} proposed a methodology for functional models with a hierarchical structure where the functions at the lowest hierarchy level are spatially correlated. \cite{DelicadoEtAl2010}, \cite{RUIZMEDINA201282}, and  \cite{Mateu2017} provided   surveys of SFD. \cite{MENAFOGLIO2017401} presented a review of complex and spatially dependent data, such as curves and surfaces. Some references to SFD with a Bayesian perspective are \cite{BaladandayuthapaniEtAl2008}, \cite{ZhangEtAl2016}, \cite{JoonEtAl2019}, and  \cite{HosseinEtAl2019}.

Functional data can also have a complex domain, e.g., a two-dimensional Euclidean domain. Spatial data can be considered 
as functional data with the same domain as the corresponding random field, for instance if observations are dense over the region or if the domain is not a subset of the Euclidean space; see, e.g.,  \cite{ALFELD19965}, \cite{Wahba81}, and \cite{Gneiting2013} for spatial data over complex domains. With the FDA approach, a realization of a random field is considered as a single point observation of the  functional data. Similarly to the case when data are curves, the continuous surface (or the manifold) needs to be estimated. For that estimation, one can use a tensor product of univariate B-splines \citep[][]{EilersMarx96,Wood2006,QingguoEtAl2010,XiaoEtAl2013}. Another way is to approximate the continuous data using finite elements analysis; see \cite{Ramsay2002}, \cite{DuchampEtAl2003}, and \cite{SangalliEtAl2013}. These types of functional data can also be observed over time; they are then called a functional time series. Functional time series can  be related to spatio-temporal data when observations are dense over the domain. Thus, the approach of functional time series can be used in complex  and large spatio-temporal data.

In this paper, our goal is to provide a review of complex and spatially correlated functional data, with two approaches using either spatial functional data or surface (manifold) data. In both approaches, we focus on covariance functions and modeling. The remainder of our paper is organized as follows: In Section \ref{SFDA}, we present basic concepts of functional random fields that include the mean, the covariance, and other important concepts related to the covariance. In this section, we also present different estimators of the various objects defined. In Section \ref{ModelingSFD}, we describe how to model small-scale and large-scale variations of spatial functional data, as well as  corresponding methods to estimate the parameters involved. Section \ref{KrigingSFD} (first approach) presents a brief description of the concept of functional kriging. In Section \ref{SurfaceD}, we describe the second approach based on surface data, which can be considered as an alternative to the analysis of spatio-temporal data, and describe how to model and estimate the continuous surface. In Section \ref{Software}, we present some software available for the analysis of spatial functional data and surface (manifold) data. In Section \ref{DataA}, we provide an example of the two approaches using a high-resolution wind speed simulated dataset in Dumat Al Jandal, Saudi Arabia. Section \ref{Discussion} ends the paper with some discussions.

\section{Functional Random Fields}\label{SFDA}
\subsection{Basic Concepts}\label{Def:FDA1}
In this section, we introduce the basic concepts of SFD. Let $(\Omega, F, P)$ be a probability space. Without loss of generality, we assume that the domain $\mathcal{T}$ of the curves is $\mathcal{T}=[0,1]$, and let $\mathcal{H}= L^{2}([0,1])$ be the Hilbert space of square integrable functions defined on $[0,1]$, equipped with the inner product $\langle f,g \rangle= \int_{0}^{1}\! f(v)g(v)\mathrm{d}v$. We denote by $\| \cdot \|_{\mathcal{H}}$ the norm in $\mathcal{H}$ induced by the inner product.  A random variable $X:\Omega \to \mathcal{H}$ taking values in the Hilbert space $\mathcal{H}$ is called a \textit{functional random variable} \citep[][]{Ramsay-Silverman2005, FerratyVieu2006, HorvathKokoszka2012}. Let $D\subset \mathbb{R}^{2}$ be a fixed  study area (e.g., a country). A random field $\{X (\mathbf{s}) : \mathbf{s}\in D\}$ taking values in $\mathcal{H}$ is called a \textit{functional random field}, that is, for each location $\mathbf{s} \in D$, $X(\mathbf{s}): \Omega \to \mathcal{H}$ is a functional random variable. 

We denote by $X(\mathbf{s}_{0}; v )$ the functional random variable at a fixed location $\mathbf{s}_{0}$ and $v\in [0,1]$, and we denote by $X(\mathbf{s}_{0}; v_{0} )$ the scalar random variable obtained by evaluating $X(\mathbf{s}; v )$ at $\mathbf{s}=\mathbf{s}_{0}$ and  $v=v_{0}$. Lastly,  we use $f$ to denote a function in $\mathcal{H}$. 

Let $X(\mathbf{s};v)$ be a functional random field with $\mathbb{E} (\|X (\mathbf{s}; v )\|_{\mathcal{H}})<\infty$, for all $\mathbf{s}\in D$. The mean $\mu(\mathbf{s}; v):= \mathbb{E} \{X(\mathbf{s}; v)\}$ of $X(\mathbf{s}; v)$ is defined  as an element of $\mathcal{H}$ such that 
\begin{equation*}
\langle \mu(\mathbf{s}; \cdot), f \rangle = \mathbb{E} \langle X(\mathbf{s}; \cdot), f \rangle, \quad \forall f\in \mathcal{H},
\end{equation*}
where the dot in $(\mathbf{s}; \cdot)$ indicates the integrated variable in $[0,1]$. This implies that $\mathbb{E}\{X(\mathbf{s}; v_{0})  \}= \mu (\mathbf{s}; v_{0})$ for almost all $v_{0}\in [0,1]$.

The covariance is one of the most studied objects in spatial statistics,  due to its relevance for prediction. If $\mathbb{E}\{\| X(\mathbf{s};v)\|^{2}_{\mathcal{H}}\}<\infty$, then the covariance operator at locations $\mathbf{s}_{1}$ and $\mathbf{s}_{2}$  is defined as an operator $C(\mathbf{s}_{1}, \mathbf{s}_{2};\sqcdot ) :\mathcal{H} \to \mathcal{H}$ such that 
\begin{align} \label{CovOp}
C(\mathbf{s}_{1}, \mathbf{s}_{2};f )(\cdot) &=  \mathbb{E}[ \langle X(\mathbf{s}_{1};\cdot) - \mu(\mathbf{s}_{1};\cdot), f \rangle  \{ X(\mathbf{s}_{2}; \cdot) -  \mu(\mathbf{s}_{2};\cdot ) \} ] \nonumber \\
	&= \int_{0}^{1} \sigma( \mathbf{s}_{1},  \mathbf{s}_{2}; u,\cdot ) f(u) \mathrm{d}u, \quad f\in \mathcal{H},
\end{align}
where $ \sigma( \mathbf{s}_{1},  \mathbf{s}_{2}; u_{0},v_{0}) := \mathbb{E} [ \{  X(\mathbf{s}_{1};u_{0}) -  \mu(\mathbf{s}_{1};u_{0} ) \}  \{X(\mathbf{s}_{2};v_{0}) -  \mu(\mathbf{s}_{2}; v_{0}) \} ]$ is the point-wise covariance, called the \textit{kernel} of $C(\mathbf{s}_{1}, \mathbf{s}_{2}; \sqcdot)$. This definition can be written in terms of a tensor operation as follows: $C(\mathbf{s}_{1}, \mathbf{s}_{2};f)= \mathbb{E} [ \{X(\mathbf{s}_{1};\cdot) - \mu(\mathbf{s}_{1};\cdot)\} \otimes \{X(\mathbf{s}_{2};\cdot) - \mu(\mathbf{s}_{2};\cdot)\} (f)  ] $.

A common assumption, in practice, is the stationarity condition of a process, which is defined as follows:
\begin{Def}[Weak stationarity] A functional random field $X(\mathbf{s};v)$ is said to be (weakly) stationary if 
\begin{enumerate}
\item $\mathbb{E}(\|X(\mathbf{s};v) \|_{\mathcal{H}}^{2} )<\infty $,
\item $\mu(\mathbf{s}; v) =  \mu(v) $, that is, the mean does not depend on the location $\mathbf{s}$, and
\item $C( \mathbf{s}_{1} + \mathbf{h}, \mathbf{s}_{2} + \mathbf{h}; \sqcdot )= C(\mathbf{s}_{1} , \mathbf{s}_{2}; \sqcdot )$, for all $\mathbf{s}_{1}, \mathbf{s}_{2}, \mathbf{h} \in D$.
\end{enumerate}
\end{Def}
The last condition is equivalent to the property that the covariance operator depends only on the increments $\mathbf{s}_{1}- \mathbf{s}_{2}$. This means, there exists a covariance operator $\tilde{C}(\mathbf{s}; \sqcdot): \mathcal{H}\to \mathcal{H}$ such that 
$$C( \mathbf{s}_{1}, \mathbf{s}_{2};f)= \tilde{C}(\mathbf{s}_{1}- \mathbf{s}_{2};f), \quad f\in \mathcal{H},$$
and so  the variance operator can be written as $\mathrm{Var}\{ X(\mathbf{s};v)\}=C(\mathbf{s}, \mathbf{s}; \sqcdot)= \tilde{C}(\mathbf{0}; \sqcdot)$. Thus, for convenience, we write $C(\mathbf{0}; \sqcdot)$ to denote the variance operator $C( \mathbf{s}, \mathbf{s};f)$ of the stationary functional random field $X(\mathbf{s}; v)$.

Now, we define the concept of isotropy for SFD.
\begin{Def}[Isotropy] A stationary functional random field $X(\mathbf{s}; v)$ is said to be isotropic if  there exists a  covariance operator $\tilde{C}_{0}(h; \sqcdot): \mathcal{H}\to \mathcal{H}$ such that  
$$C( \mathbf{s}_{1}, \mathbf{s}_{2}; f )= \tilde{C}_{0}(h;f), \quad f\in \mathcal{H},$$ 
where $h= \| \mathbf{s}_{1}- \mathbf{s}_{2}\|$, for all $\mathbf{s}_{1}, \mathbf{s}_{2} \in D$.
\end{Def}

In spatial statistics, the variogram plays an important role to make inference. The extension of this concept to functional random fields is as follows: The variogram operator $\Gamma$ is defined as the variance operator of the difference between the functional random field at two locations $ \mathbf{s}_{1}, \mathbf{s}_{2}$, i.e., $\Gamma( \mathbf{s}_{1}, \mathbf{s}_{2}; \sqcdot ) :=\frac{1}{2} \mathrm{Var} \{X(\mathbf{s}_{1};v) - X(\mathbf{s}_{2};v) \}$. If the functional random field $X(\mathbf{s};v)$ has a finite second moment, then we have  
$
2 \Gamma( \mathbf{s}_{1}, \mathbf{s}_{2}; f) = C(\mathbf{s}_{1}, \mathbf{s}_{1}; f) +
C(\mathbf{s}_{2}, \mathbf{s}_{2};f) - C(\mathbf{s}_{1}, \mathbf{s}_{2};f) - C(\mathbf{s}_{2},\mathbf{s}_{1};f),
$ 
for $f\in \mathcal{H}$. Thus, if the functional random field is stationary, then there exists a variogram operator $\tilde{\Gamma} (\mathbf{s}; \sqcdot):\mathcal{H}\to \mathcal{H}$ such that it satisfies
\begin{equation}\label{VarOp}
\Gamma( \mathbf{s}_{1}, \mathbf{s}_{2};f)= \tilde{\Gamma}( \mathbf{s}_{1}- \mathbf{s}_{2};f)=  C(\mathbf{0};f)  -  C(\mathbf{s}_{1}- \mathbf{s}_{2};f), \quad f\in \mathcal{H},
\end{equation}
where the corresponding kernel is $\gamma( \mathbf{s}_{1}- \mathbf{s}_{2}; u,v) =  \sigma (\mathbf{0}; u,v ) - \sigma( \mathbf{s}_{1}- \mathbf{s}_{2}; u,v) $. Furthermore, if $X(\mathbf{s};v)$ is isotropic, then there exists $\tilde{\Gamma}_{0}(h; \sqcdot): \mathcal{H}\to \mathcal{H}$ such that 
$$\Gamma( \mathbf{s}_{1}, \mathbf{s}_{2};f) = \tilde{\Gamma}_{0}( h;f), \quad f\in \mathcal{H},$$
where $h= \| \mathbf{s}_{1}- \mathbf{s}_{2}\|$, for all $\mathbf{s}_{1}, \mathbf{s}_{2} \in D$.

\begin{rem}  Unlike the finite dimensional case (multivariate random field) where the covariance is composed of matrices, here,  the covariance is composed  of operators, since the space of functional data is intrinsically infinite-dimensional.
\end{rem}


Another way to describe the second-order spatial dependence of the functional random field is by using a ``global'' measure. This global measure is the trace-covariogram  $\sigma_{\mathrm{tr}} : D\times D\to \mathbb{R}$  \citep[][]{GiraldoEtAl2011,menafoglio2013} defined as
\begin{equation}\label{tCov}
 \sigma_{\mathrm{tr}} ( \mathbf{s}_{1}, \mathbf{s}_{2})= \mathbb{E} \{ \langle X(\mathbf{s}_{1}; \cdot) - \mu (\mathbf{s}_{1}; \cdot), X(\mathbf{s}_{2};\cdot) -  \mu (\mathbf{s}_{2}; \cdot) \rangle \}=\int_{0}^{1}\! \sigma (\mathbf{s}_{1}, \mathbf{s}_{2}; v,v ) \mathrm{d}v.
\end{equation}
The trace-covariogram computes the covariance of the inner product of the functional random field at two locations. Thus, it summarizes the covariance on the diagonal, and so, in general, it depends only on the locations. If the functional random field is stationary, then there exists $\tilde{\sigma}_{\mathrm{tr}}: D \to \mathbb{R}$ such that it depends only on the separation vector $ \mathbf{s}_{1}- \mathbf{s}_{2}$, that is, 
$$\sigma_{\mathrm{tr}} ( \mathbf{s}_{1}, \mathbf{s}_{2}) = \tilde{\sigma}_{\mathrm{tr}} (  \mathbf{s}_{1}- \mathbf{s}_{2}).$$
In addition, if the functional random field is isotropic, then there exists $\tilde{\sigma}_{\mathrm{tr},0}: \mathbb{R}\to \mathbb{R}$ such that 
$$ \sigma_{\mathrm{tr}} ( \mathbf{s}_{1}, \mathbf{s}_{2}) = \tilde{\sigma}_{\mathrm{tr},0} (h),$$
where $h= \| \mathbf{s}_{1}- \mathbf{s}_{2}\|$.

Similarly, the trace-variogram is defined in terms of the inner product of the difference, i.e.,
\begin{equation*}
 \gamma_{\mathrm{tr}} (  \mathbf{s}_{1}- \mathbf{s}_{2}) =  \frac{1}{2} \mathbb{E} \left \{ \langle X(\mathbf{s}_{1}; \cdot) - X(\mathbf{s}_{2};\cdot), X(\mathbf{s}_{1};\cdot ) -  X(\mathbf{s}_{2}; \cdot)\rangle   \right\}  \\
     -   \frac{1}{2}\| \mu (\mathbf{s}_{1};\cdot)-  \mu (\mathbf{s}_{2};\cdot) \|_{\mathcal{H}}^{2}.
\end{equation*}
We observe that, if  $X(\mathbf{s}; v)$ is stationary, then   $\sigma (\mathbf{s}; u,v) = \sigma (\mathbf{0}; u,v) - \gamma (\mathbf{s}; u,v).$ Thus, the  trace-variogram satisfies
$$\gamma_{\mathrm{tr}} (  \mathbf{s}_{1}- \mathbf{s}_{2})  =  \sigma_{\mathrm{tr}} (\mathbf{0} )- \sigma_{\mathrm{tr}}( \mathbf{s}_{1}- \mathbf{s}_{2} ) .$$

Trace-covariogram and trace-variogram are also important for optimization problems. Especially if we want to use the criterion of minimizing  equations of the form  $\mathbb{E}(\langle X, Y \rangle)$, as in \eqref{OptProblem2} below.

\subsection{Estimation}\label{sec:Est}
Now, we describe estimators  of the mean $\mu$, the covariance $C$, and the variogram $\Gamma$. For this purpose, we assume that $X(\mathbf{s};v)$ is an isotropic (stationary) functional random field with mean $\mu(v)$ and covariance operator $C$.  Let $x(\mathbf{s}_{1}; v), \ldots, x(\mathbf{s}_{n}; v)$ be observations of the functional random field $X(\mathbf{s};v)$.   We assume that the observations $x(\mathbf{s}_{i}; v)$ are given in the functional form. Although in real data, $x(\mathbf{s}_{i}; v)$ are observed on a finite set of points $v_{i1},\ldots, v_{im}$, the continuous curves should be estimated \citep[see][]{Ramsay-Silverman2005}. 

The main feature of spatial data is that ``nearby'' data look similar, and an estimator must take into account such spatial dependence.  Otherwise,  it will not have desirable properties, such as consistency. 

\subsubsection{Mean estimation}

We describe two different approaches to obtain an estimator of the mean $\mu(v)$ \citep{GromenkoAndKokoszka2012A}. A model of the mean can be written as 
\begin{equation}\label{Model-mu}
X(\mathbf{s}; v) = \mu(v) + \varepsilon(\mathbf{s}; v),
\end{equation}
where $ \varepsilon(\mathbf{s}; v)$ is an isotropic functional random field with zero mean and covariance operator $C$. 

The first approach is similar to the kriging method. Specifically, this is defined as a weighting of the observed curves:
\begin{equation}\label{MeanEst}
\hat{\mu}(v)= \sum_{i=1}^{n} w_{i}x(\mathbf{s}_{i}; v) ,
\end{equation}
where the weights $w_{i}$ are estimated by solving the optimization problem
\begin{equation}\label{OptProblem}
\min_{w_{1}, \ldots, w_{n} } \mathbb{E}\{\langle  \hat{\mu}-\mu , \hat{\mu}-\mu  \rangle \}= \min_{w_{1}, \ldots, w_{n} } \mathbb{E} \left\{ \left \| \sum_{i=1}^{n} w_{i}x(\mathbf{s}_{i}; \cdot)-\mu \right \|^{2}_{\mathcal{H}} \right \} ,
\end{equation}
subject to the condition $\sum_{i=1}^{n}  w_{i} =1$. Using the Lagrange multiplier method, this leads to solve 
\begin{equation}\label{OptProblem2}
\sum_{i=1}^{n} w_{i}=1, \quad \sum_{i=1}^{n} w_{i}\,  \sigma_{\mathrm{tr},\varepsilon} (\mathbf{s}_{i}, \mathbf{s}_{j} )-\lambda=0, \quad j=1,\ldots, n ,
\end{equation}
where $ \sigma_{\mathrm{tr},\varepsilon}(\mathbf{s}_{1}, \mathbf{s}_{2} )$ is the trace-covariogram of $\varepsilon(\mathbf{s};v)$. Thus, the  estimation problem \eqref{OptProblem} becomes estimating the matrix $\{  \sigma_{\mathrm{tr},\varepsilon} (\mathbf{s}_{i}, \mathbf{s}_{j} )\}_{i,j=1}^{n}$. Since the functional random field $ \varepsilon(\mathbf{s}; v)$ is unobserved, a common approach is to use an iterative procedure.  At the first iteration, an initial estimator of $\mu(v)$ is obtained by assuming that $\varepsilon(\mathbf{s}; v)$ is spatially uncorrelated, i.e., $\hat{\mu}_{0}(v)= \frac{1}{n} \sum_{i=1}^{n} x(\mathbf{s}_{i}; v)$. Next, $\hat{\mu}_{0}(v)$ is subtracted from the data $x(\mathbf{s}_{i}; v)$, then an initial estimator of $ \sigma_{\mathrm{tr},\varepsilon} (\mathbf{s}_{i}, \mathbf{s}_{j} )$ is obtained as described below in \eqref{trace-cov}. At the second iteration, the mean is re-estimated by solving \eqref{OptProblem2}, with the initial information of  $ \sigma_{\mathrm{tr},\varepsilon} (\mathbf{s}_{i}, \mathbf{s}_{j} )$. This process is repeated until convergence. 
 
The second approach uses finite basis functions similarly to the cokriging method on the coefficients of the basis functions, see \cite{Goulard1993}, \cite{NERINI2010}, and \cite{GiraldoEtAl2011}. Let $ \{\eta_{1}(v), \ldots, \eta_{K}(v) \}$ be  basis functions, e.g., Fourier basis functions or B-spline basis functions.  Then, the observed curves are approximated  as
\begin{equation}\label{approxX}
x(\mathbf{s}_{i}; v) \approx \sum_{k=1}^{K} z_{k} (\mathbf{s}_{i}) \eta_{k}(v) ,\quad i=1, \ldots, n,
\end{equation}
where  $z_{k} (\mathbf{s}_{i}) = \langle  x(\mathbf{s}_{i}; \cdot), \eta_{k} \rangle$, which is a scalar for each $i$ and $k$. Let $Z_{k}(\mathbf{s})$ be the corresponding scalar random field with realization $\{z_{k}(\mathbf{s}_{i})\}_{i=1}^{n}$. Then, by using \eqref{Model-mu}, the mean can be approximated as
$$\mu(v)\approx \sum_{k=1}^{K} \mathbb{E}\{ Z_{k} (\mathbf{s})\} \eta_{k}(v),$$
and $\mathbb{E}\{ Z_{k} (\mathbf{s})\}$ should be estimated. For this,  from  \eqref{Model-mu} we have that
$$Z_{k} (\mathbf{s}_{i}) = \mu^{z}_{k} + \varepsilon^{k}(\mathbf{s}_{i}), \quad i=1, \ldots, n,$$ 
where $ \mu_{k}^{z}  = \langle \mu,\eta_{k} \rangle= \mathbb{E}\{ Z_{k} (\mathbf{s})\}$, and  $\varepsilon^{k}(\mathbf{s}_{i})= \langle  \varepsilon(\mathbf{s}_{i}; \cdot), \eta_{k}\rangle$, for each $k=1,\ldots,K$ .
Thus, $Z_{k} (\mathbf{s}_{i})$ is an isotropic (stationary) scalar random field, for $k=1,\ldots,K$. 
Once $\mu^{z}_{k}$ is estimated, we obtain the estimator of the mean 
$$\hat{\mu}(v) = \sum_{k=1}^{K} \hat{\mu}^{z}_{k}\eta_{k}(v).$$
To estimate $\mu^{z}_{k}$, we can use the ordinary least squares (OLS) estimator, by minimizing
$$\sum_{i=1}^{n} \left\{ z_{k}(\mathbf{s}_{i})- w_{i}\mu^{z}_{k} \right\}^{2},$$
with $w_{i}=1$ for all $i=1, \ldots, n$, or a weighted least squares (WLS) estimator, i.e, $w_{i}$ are estimated with the information of the  covariance matrix of $\varepsilon^{k}(\mathbf{s})$. Also, we can use a generalized least squares (GLS) estimator, which can perform better than OLS or WLS; see \cite{Schabenberger2005} and  \cite{Cressie2015} for the development of these estimators.

 Another way to estimate $\mu^{z}_{k}$ is by using Maximum  Likelihood Estimators (MLEs). For this, assume that, for each $i$, $X(\mathbf{s}_{i}; v)$ is a Gaussian process in $[0,1]$. Then, for each $k$, $\mathbf{z}_{k}= \{  z_{k} (\mathbf{s}_{1} ), \ldots, z_{k} (\mathbf{s}_{n}) \}^{\mathsf{T}}$ is a realization of 
$\mathbf{Z}_{k} \sim {\cal N}_{n} (\mu^{z}_{k}\mathbf{1} , \boldsymbol{\Sigma}(\boldsymbol{\theta})),$ 
where $\boldsymbol{\Sigma}(\boldsymbol{\theta})$ is the covariance of $\varepsilon^{k}(\mathbf{s})$ described by some valid parametric covariance model, e.g., Mat\'ern, exponential, or spherical. Then, the log likelihood is 
\begin{equation}\label{LogL-1}
l (\boldsymbol{\theta}, \mu^{z}_{k} )= - \frac{n}{2} \log (2\pi) - \frac{n}{2} \log\, \mathrm{det} \{\boldsymbol{\Sigma}(\boldsymbol{\theta}) \} - \frac{1}{2} (\mathbf{z}_{k} - \mu^{z}_{k} \mathbf{1}  )^{\mathsf{T}} \boldsymbol{\Sigma}(\boldsymbol{\theta})^{-1} (\mathbf{z}_{k} - \mu^{z}_{k} \mathbf{1}  ).
\end{equation}
To read more details of this estimator, see \cite{Stein99}.

The coefficients in estimator \eqref{MeanEst} can be also considered as operators $\omega_{i}: \mathcal{H}\to \mathcal{H}$. In this case the estimator of the mean takes the form $\hat{\mu}(v)= \sum_{i=1}^{n} \omega_{i}\{ x(\mathbf{s}_{i}; \cdot) \}(v)$, where the coefficients, which are  operators, need to be estimated under a certain constraint on $\sum_{i=1}^{n} \omega_{i}$, see \cite{NERINI2010}. This approach takes into account the information of the whole curve to define a specific weight for each point $v\in [0,1]$. Thus, it is expected to obtain better results than estimator \eqref{MeanEst}, but the estimation procedure can be complicated.
	
In the context of a nonstationary functional random field, the mean can depend on the location $\mathbf{s}$, in this case,  $\mu(\mathbf{s}; v)$ can be represented by a linear model $\mu(\mathbf{s};v)= \sum_{l=0}^{L} a_{l}( \mathbf{s}) f_{l}(v)$, where $f_{l}(v)$, $l=0,\ldots,L$, are elements of $\mathcal{H}$ independent of the spatial location $\mathbf{s}$, and $a_{0}(\mathbf{s})=1$ \citep[][]{menafoglio2013,Caballero2013}. The latter approach  is described with more details in Section~\ref{section:large-v}.

\subsubsection{Covariance estimation}\label{EstCov}

Now, we describe how to estimate the covariance. We assume that $\mathbb{E}\{ X( \mathbf{s}; v)\}=0$. Let $N(h)$ be a set of pairs of indexes defined as $N(h)= \{ (i,j) : \| \mathbf{s}_{i}-\mathbf{s}_{j}\|= h  \}$, and let $| N(h)|$ be the cardinality of $N(h)$. The empirical covariance operator of the functional random field is defined as 
\begin{equation*}
\hat{C}(h;f) = \frac{1}{| N(h)|} \sum_{(i,j) \in N(h)} x(\mathbf{s}_{i}) \otimes x(\mathbf{s}_{j}) (f).
\end{equation*}
In practice, the distance between $ \mathbf{s}_{i}$ and $\mathbf{s}_{j}$ is not considered to be exactly $h$, instead $ \| \mathbf{s}_{i}-\mathbf{s}_{j}\| \in (h-\delta, h+\delta) $, with $\delta >0$. Also, since it is almost impossible to obtain $\hat{C}(h; f) $ for all $h$, discretized values of $h$, $h_{1}, \ldots, h_{m}$, are computed.

The usual approach in the scalar (multivariate) random field case is to fit a valid parametric model to the empirical covariance by least squares methods \citep[][]{Cressie2015}. However, in spatial functional data, for each $h$, $\hat{C}(h;\sqcdot)$ is an operator. Thus, it requires new mathematical developments to define a ``covariance model'' in this context.  Finite basis functions have been used to overcome the modeling of the covariance function in $\mathcal{H}$ \citep[][]{NERINI2010}. Recall that $C(h; \sqcdot)$ is represented by its kernel $\sigma(h; u,v)$ (see \eqref{CovOp}). Now, we assume that $x(\mathbf{s}_{i}; v) $ is approximated by a finite set of basis functions, as in \eqref{approxX}. Then, we obtain 
$$\sigma(h; u,v)= \mathbb{E} \{ X(\mathbf{s+h}) X(\mathbf{s}) \}= \boldsymbol{\eta}^{\mathsf{T}}(u) \mathbb{E}\{\mathbf{Z}(\mathbf{s+h})\mathbf{Z}^{\mathsf{T}}(\mathbf{s}) \}  \boldsymbol{\eta}(v),$$
where $\boldsymbol{\eta}(v)= \{\eta_{1}(v), \ldots, \eta_{K}(v) \}^{\mathsf{T}}$, $\mathbf{Z}(\mathbf{s})=\{Z_{1}(\mathbf{s}), \ldots, Z_{K}(\mathbf{s})\}^{\mathsf{T}} $, and $Z_{k}(\mathbf{s})= \langle X(\mathbf{s}; \cdot), \eta_{k} \rangle$ is the  scalar random field with mean zero, and so  $\mathbf{Z} (\mathbf{s})$  is a multivariate random field with covariance 
$$\boldsymbol{\Sigma}(h)= \mathrm{Cov} \{\mathbf{Z}(\mathbf{s+h}), \mathbf{Z}(\mathbf{s}) \}= \{ \boldsymbol{\Sigma}_{kl}(h) \}_{k,l=1}^{K},$$
where $\boldsymbol{\Sigma}_{kl}(h)=  \mathrm{Cov} \{Z_{k}(\mathbf{s+h}), Z_{l}(\mathbf{s}) \}$ for $k,l=1, \ldots, K$, $h=\|\mathbf{h}\|$. Thus, estimating the covariance operator $C$ is equivalent to estimating the covariance of a multivariate random field $\mathbf{Z} (\mathbf{s})$, which can be modeled with a valid covariance function $\boldsymbol{\Sigma}(h;\boldsymbol{\theta})$ \citep[][]{GentonAndKleiber2015}.

If $\boldsymbol{\eta}(v)$ is a set of orthogonal basis functions, then we can estimate each marginal-covariance functions $\boldsymbol{\Sigma}_{kk}(h)$, for $k=1,\ldots, K$, separately. In addition, if $X(\mathbf{s};v)$ is assumed to be a Gaussian process, we can use the log likelihood \eqref{LogL-1} with $\mu^{z}_{k}=0$, i.e., maximizing $l (\boldsymbol{\theta}, 0 )$ as function of $\boldsymbol{\theta}$. Also, the restricted maximum likelihood (REML) estimation can be used in this case, see \cite{Stein99}.

Once the estimator $\hat{\boldsymbol{\Sigma}}(h;\boldsymbol{\theta})$ is obtained, then we define an estimator of $\sigma(h; u,v)$ as follows,
\begin{equation}\label{kernelEst}
\hat{\sigma}(h; u,v) = \boldsymbol{\eta}^{\mathsf{T}}(u) \hat{\boldsymbol{\Sigma}}(h;\boldsymbol{\theta})  \boldsymbol{\eta}(v)
\end{equation}
and then, $\hat{C}(h;f)(\cdot)= \int_{0}^{1} \hat{\sigma}(h; u,\cdot) f(u)\mathrm{d}u$.

In the estimation of the covariance of $Z_{k}$, we can also use a nonparametric approach \citep{HallEtAl94,Hall1994}. Bayesian approaches  can be found in, e.g.,  \cite{BanerjeeEtAl2015} and \cite{DiggleEtAl2007}.

As mentioned before, the trace-covariogram is a measure that describes  dependence globally, in the sense that it integrates the kernel $\sigma(\mathbf{s}; u,v)$ on the diagonal. Also, the trace-covariogram appears in several optimization problems \citep[][]{DelicadoEtAl2010, GiraldoEtAl2011, menafoglio2013}, that make its estimation important. The empirical trace-covariogram is defined as 
\begin{equation}\label{trace-cov}
 \hat{\sigma}_{\mathrm{tr}}(h) = \frac{1}{ | N(h) |} \sum_{ (\mathbf{s}_{i},\mathbf{s}_{j}) \in N(h) } \int_{0}^{1}  x(\mathbf{s}_{i}; v)  x(\mathbf{s}_{j}; v) \mathrm{d}v.
\end{equation}
This empirical estimator is computed on discrete values $h_{1}, \ldots, h_{m}$ of $h$. Then, we can fit any covariance model by least squares methods to the empirical estimates $ \hat{\sigma}_{\mathrm{tr}}(h_{1}), \ldots  , \hat{\sigma}_{\mathrm{tr}}(h_{m})$.

We observe that estimating the trace-covariogram is much easier than estimating the covariance $C$, where the inner products $\langle  x(\mathbf{s}_{i}; \cdot), x(\mathbf{s}_{j}; \cdot)  \rangle=\int_{0}^{1}  x(\mathbf{s}_{i}; v)  x(\mathbf{s}_{j}; v) \mathrm{d}v $ can be computed using the \texttt{R} \citep[][]{RSoft} package \textit{fda} \citep[][]{fda}.

\subsubsection{Variogram estimation}

One of the advantages of using the variogram is the  robustness under misspecification of the mean. Also, the variogram can be used  to estimate the covariance. The empirical variogram operator of the functional random field is defined as
\begin{equation*}
\hat{\Gamma}(h;f) = \frac{1}{2 | N(h)|} \sum_{(i,j) \in N(h)} \{ x(\mathbf{s}_{i})-x(\mathbf{s}_{j})\} \otimes \{ x(\mathbf{s}_{i})-x(\mathbf{s}_{j})\} (f), \quad f\in \mathcal{H}.
\end{equation*}
Similarly to the estimator of the covariance operator, if a basis function $\boldsymbol{\eta}(v)$ is assumed, then the corresponding kernel $\gamma$ is obtained as  
\begin{align*}
2 \gamma(h; u,v) &=  \mathbb{E} [ \{ X(\mathbf{s+h})- X(\mathbf{s}) \}^{2} ] \\
                         &=\boldsymbol{\eta}^{\mathsf{T}}(u) \mathbb{E} [ \{\mathbf{Z}(\mathbf{s}+\mathbf{h})-\mathbf{Z}(\mathbf{s}) \} \{\mathbf{Z}(\mathbf{s}+\mathbf{h})-\mathbf{Z}(\mathbf{s}) \}^{\mathsf{T}}]  \boldsymbol{\eta}(v),
\end{align*}
where $h= \|\mathbf{h}\|$. Thus, we need to estimate the variogram of the multivariate random field $\mathbf{Z}(\mathbf{s}) = \{Z_{1}(\mathbf{s}), \ldots, Z_{K}(\mathbf{s}) \}^{\mathsf{T}} $, to obtain the estimator of the kernel $\gamma(h;u,v)$.

The corresponding  empirical trace-variogram is defined as
\begin{equation}\label{EmTraVar}
 \hat{\gamma}_{\mathrm{tr}}(h)= \frac{1}{2 | N(h) |} \sum_{ (i,j) \in N(h) } \int_{0}^{1}  \{ x(\mathbf{s}_{i}; v) -  x(\mathbf{s}_{j}; v) \}^{2} \mathrm{d}v.
\end{equation}
Then, following the common method used in spatial statistics, a variogram model is fitted to $\hat{\gamma}_{\mathrm{tr}}(h_{1}),\ldots, \hat{\gamma}_{\mathrm{tr}}(h_{m})$.

The principal component analysis is generally important in statistics because of its applicability to dimensional reduction techniques. In the context of functional data, a function $\zeta \in \mathcal{H}$ is called \textit{eigenfunction} of the operator $C(\mathbf{s},\mathbf{s};\sqcdot)$, if $C(\mathbf{s}, \mathbf{s}; \zeta )=\lambda \zeta $ with $\lambda$ a positive real number. The estimation of eigenfunctions can be obtained from the estimated covariance operator $\hat{C}(\mathbf{s}, \mathbf{s}; \sqcdot)$. If we are only interested in the eigenfunctions,  we can obtain the estimators  without estimating the covariance operator, which allows us to reduce computational costs \citep[][]{ZhouEtAl2010, GromenkoAndKokoszka2012A}. Also see \cite{LiuEtAl2017} for functional principal component analysis of spatial functional data.

\section{Modeling Functional Random Fields}\label{ModelingSFD}

Spatial statistics studies the variations among the observed data at different locations. The spatial variation is generally described through the mean and the covariance \citep{haining2003}. The mean represents the large-scale variations, and the covariance represents the small-scale variations. 

In this section, we describe the statistical models for spatial functional data. We denote the observed functional data as  $y(\mathbf{s}_{1};v), \ldots, y(\mathbf{s}_{n};v)$, and we denote by  $Y(\mathbf{s};v)$ the corresponding functional random field with realization $y(\mathbf{s}_{i};v)$, $i=1,\ldots, n$. The functional random field  $X(\mathbf{s};v)$ denotes a stationary functional random field with covariance function $C_{X}$, and $\varepsilon(\mathbf{s};v)$ denotes a functional white noise, i.e.,  $\mathbb{E}\{\varepsilon(\mathbf{s},v)\}=0$, with a covariance such that $C_{\varepsilon} (\mathbf{s}_{1}, \mathbf{s}_{2}; \sqcdot )=0$, if $\mathbf{s}_{1}\neq \mathbf{s}_{2}$.

\subsection{Large-scale Variation}\label{section:large-v}

Regression models with covariates in terms of the locations can be used to model large-scale spatial variations. Such  regression models must account for spatial dependence. In general, estimators obtained from regression models are smooth functions defined in $D$. In the context of spatial functional data, functional regression models in $\mathcal{H}$ extend models of finite-dimensional data to model large-scale variations. In \cite{Caballero2013}, \cite{menafoglio2013}, and \cite{ReyesEtAl2015}, covariates are assumed  to be separable in the spatial component and the continuity of the data. In this case, the model for large-scale variations is 
\begin{equation}\label{LD-model}
Y(\mathbf{s};v) = \sum_{l=0}^{L} a_{l}( \mathbf{s}) f_{l}(v) + \varepsilon(\mathbf{s};v), 
\end{equation}
where $f_{l}\in \mathcal{H}$ are independent of $\mathbf{s}$, $a_{0}(\mathbf{s}):=1$, $\{a_{l}(\mathbf{s})\}_{l=1}^{L}$ are known scalar regressors, and  $\varepsilon(\mathbf{s}; v)$ is the functional white noise. For example, with $L=5$, we could specify the scalar regressors $a_{l}(\mathbf{s})$ as
$$
a_{1}(\mathbf{s})=s_{1},\, a_{2}(\mathbf{s})=s_{2},\, a_{3}(\mathbf{s})= s_{1}s_{2}, \, a_{4}(\mathbf{s})= s_{1}^{2}, \mbox{ and } a_{5}(\mathbf{s})=s_{2}^{2},
$$ 
where $\mathbf{s}=(s_{1}, s_{2})$ denotes the coordinates of a spatial location. In this example, $f_{0}(v), \ldots, f_{5}(v)$ need to be estimated. Thus, model \eqref{LD-model} is a functional regression model with scalar covariates $a_{l}(\mathbf{s})$, $l=0,1,\ldots,L$.

The mean of the functional random field $Y(\mathbf{s};v)$ in \eqref{LD-model}  is
$$
\mathbb{E}\{Y(\mathbf{s};v)\}= \mu(\mathbf{s};v) = \sum_{l=0}^{L} a_{l}( \mathbf{s}) f_{l}(v),
$$
and the covariance $C_{Y}$ of  $Y(\mathbf{s};v) $ is $C_{Y}= C_{\varepsilon}$, which is zero at $(\mathbf{s}_{i}, \mathbf{s}_{j})$ if $i \neq j$. Then, the spatial variation is described through the mean of the functional random field $Y(\mathbf{s};v)$. The covariates  $a_{l}(\mathbf{s})$ capture the spatial dependence, and the set of functions $\{f_{l}(v)\}_{l=0}^{L}$ carries the continuity of the functional data.

The estimators of  $f_{l}(v)$ in model \eqref{LD-model} can be obtained using the OLS method. The matrix form of the model can then be written as
$$
\mathbf{y}(v)= \mathbf{Af}(v) + \boldsymbol{\varepsilon}(v),
$$
where $\mathbf{y}(v)= \{ y(\mathbf{s}_{1};v), \ldots, y(\mathbf{s}_{n};v)   \}^{\mathsf{T}}$, $\mathbf{A}= \{ a_{l}(\mathbf{s}_{i}) \}_{i,l}$, $i=1,\ldots, n$, $l=0,1,\ldots, L$,  is the design matrix, $\mathbf{f}(v)= \{f_{0}(v), \ldots, f_{L}(v) \}^{\mathsf{T}}$, and  $\boldsymbol{\varepsilon}(v)= \{\varepsilon(\mathbf{s}_{1};v), \dots, \varepsilon(\mathbf{s}_{n};v) \}^{ \mathsf{T}  }$.  Then, the OLS estimator is obtained by solving the optimization problem 
\begin{equation}\label{Op-LD}
\min_{f_{0}, \ldots, f_{L}\in \mathcal{H}} \sum_{i=1}^{n}\left \| y(\mathbf{s}_{i};\cdot) -  \sum_{l=0}^{L} a_{l}( \mathbf{s}_{i}) f_{l} \right \|^{2}_{\mathcal{H}}.
\end{equation}
Under some conditions \citep[][]{menafoglio2013}, \eqref{Op-LD} admits a unique solution 
 $$\hat{\mathbf{f}}(v)= ( \mathbf{A}^{\mathsf{T}}\mathbf{A})^{-1} \mathbf{A}^{\mathsf{T}} \mathbf{y}(v).$$ 
Thus, the drift estimator is obtained as 
 \begin{equation}\label{estimatorLD}
 \hat{\boldsymbol{\mu}}(v) = \mathbf{A} (\mathbf{A}^{\mathsf{T}}\mathbf{A})^{-1} \mathbf{A}^{\mathsf{T}} \mathbf{y}(v).
 \end{equation}
Since the estimator $ \hat{\boldsymbol{\mu}}(v) $ in \eqref{estimatorLD} is a linear combination of the observed curves, $ \hat{\boldsymbol{\mu}}(v) $  inherits the continuity property of $\mathbf{y}(v)$. 

The basis functions approach offers another alternative to obtain estimators of $\{f_{l}(v)\}$ \citep[][]{ReyesEtAl2015}. Each component of the model \eqref{LD-model} can be assumed to be in the space generated by finite basis functions,  i.e., 
\begin{equation*}
y(\mathbf{s}_{i};v)=\sum_{k=1}^{K}z_{ik} \eta_{k}(v),\,  f_{l}(v)= \sum_{k=1}^{K} b_{lk}\eta_{k}(v), \mbox{ and } \varepsilon(\mathbf{s}_{i};v)=\sum_{k=1}^{K}e_{ik} \eta_{k}(v),
\end{equation*}
for $i=1, \ldots,n$, and $l=0,1, \ldots, L$, where $\boldsymbol{\eta}(v)= \{\eta_{1}(v), \ldots, \eta_{K}(v) \}^{\mathsf{T}}$ is the basis function.  In this case, the matrix form of the model \eqref{LD-model} is 
\begin{equation*}
\mathbf{Z} \boldsymbol{\eta}(v)= \mathbf{AB}\boldsymbol{\eta}(v) + \mathbf{E} \boldsymbol{\eta}(v),
\end{equation*}
where $\mathbf{Z}=\{z_{ik}\}$, $\mathbf{B}=\{b_{lk}\}$, and $\mathbf{E}=\{e_{ik}\}$, $i=1, \ldots, n$, $l=0,1, \ldots, L$, and $k=1,\ldots, K$. The corresponding normal equation is 
$$\mathbf{A}^{\mathsf{T}} \mathbf{Z} \mathbf{J}_{\boldsymbol{\eta}}= \mathbf{A}^{\mathsf{T}} \mathbf{A B} \mathbf{J}_{\boldsymbol{\eta}},$$
where $ \mathbf{J}_{\boldsymbol{\eta}}= \int \boldsymbol{\eta}(v)\boldsymbol{\eta}^{\mathsf{T}}(v)\mathrm{d}v$. The solution for $\mathbf{B}$ is found by vectorizing the normal equation, that is
\begin{equation}\label{estimatorLD2}
\mathrm{vec}(\hat{\mathbf{B}})= ( \mathbf{J}_{\boldsymbol{\eta}}^{\mathsf{T}}\otimes \mathbf{A}^{\mathsf{T}}\mathbf{A} )^{-1}  \mathrm{vec}(\mathbf{A}^{\mathsf{T}} \mathbf{Z} \mathbf{J}_{\boldsymbol{\eta}}).
\end{equation}
Consequently, $ \hat{\boldsymbol{\mu}}(v) = \mathbf{A}\hat{\mathbf{B}}\boldsymbol{\eta}(v)$. 

In the context of explicative modeling, we can use other spatial functional data as covariates to describe the mean $\mu(\mathbf{s}; v)$ \citep[][]{Ignaccolo2014}. For example, the mean can be modeled as $\mu(\mathbf{s};v)= \beta_{0}(v)+ \sum_{p=1}^{P} \beta_{p}(v) U_{p}(\mathbf{s};v)$, where $\{U_{p}(\mathbf{s};v)\}_{p=1}^{P}$ are the functional covariates, and $ \beta_{0}(v),  \beta_{1}(v), \ldots,  \beta_{P}(v)$ are the functional  parameters to be estimated.

\subsection{Small-scale Variation}\label{section:small-v}

Small-scale variations are usually represented through the covariance structure. The modeling of covariance is one of the most studied subjects in spatial statistics \citep[][]{Stein99, Cressie2015, GentonAndKleiber2015}. For finite-dimensional data, the exponential, the Gaussian, and the Mat\'ern are examples of parametric covariance function  models.  

In this section, we assume that the mean is constant over locations, and without loss of generality set to be zero. A model to describe the small-scale variations is 
\begin{equation}\label{small-v-m}
Y(\mathbf{s};v) = X(\mathbf{s};v) + \varepsilon(\mathbf{s};v),
\end{equation}
where $\varepsilon(\mathbf{s};v)$ represents the functional white noise, and is assumed to be uncorrelated with $X(\mathbf{s};v)$.

The mean of $Y(\mathbf{s};v)$ in model \eqref{small-v-m} is zero, and the covariance  is such that
\begin{equation}\label{Mod:CovOp}
C_{Y}(\mathbf{s}_{i}, \mathbf{s}_{j}; f) = C_{X}(\mathbf{s}_{i}, \mathbf{s}_{j};f)+ \mathds{1} (i=j) C_{\varepsilon}(\mathbf{s}_{i}, \mathbf{s}_{j};f), \quad f\in \mathcal{H},
\end{equation}
where $\mathds{1}(\cdot)$ is the indicator function. Thus, we need to consider the additional term $C_{\varepsilon}(\mathbf{s}, \mathbf{s}; \sqcdot)$ in the estimation. Similarly, as in Section \ref{sec:Est}, a basis function approach can be used to obtain the estimators. 

Let   $\mathbf{Z}(\mathbf{s})=\{Z_{1}(\mathbf{s}), \ldots, Z_{K}(\mathbf{s})\}^{\mathsf{T}}$ be the multivariate random field obtained from the projection of $Y(\mathbf{s};v)$ onto the basis functions $\{\eta_{k}\}_{k=1}^{K}$, that is, each component of $\mathbf{Z}(\mathbf{s})$ is defined as $Z_{k}(\mathbf{s})= \langle Y(\mathbf{s}; \cdot),\eta_{k} \rangle$, which are scalar random fields with mean zero. From \eqref{small-v-m}, for each $k=1,\ldots,K$, we have that the process $Z_{k}(\mathbf{s})$ is such that
$$Z_{k} (\mathbf{s})= \langle X(\mathbf{s}), \eta_{k} \rangle +  \langle \varepsilon(\mathbf{s}), \eta_{k} \rangle.$$
The variance of $Z_{k}(\mathbf{s})$ is 
\begin{align*}
\mathbb{E} \{  \langle Y(\mathbf{s}; \cdot),\eta_{k} \rangle  \langle Y(\mathbf{s}; \cdot),\eta_{k} \rangle \} &=  \langle C_{Y}(\mathbf{s},\mathbf{s}; \eta_{k}), \eta_{k} \rangle \\
  & = \langle C_{X}(\mathbf{s},\mathbf{s}; \eta_{k}), \eta_{k} \rangle + \langle C_{\varepsilon}(\mathbf{s},\mathbf{s}; \eta_{k}), \eta_{k} \rangle,
\end{align*}
where the last equality is obtained using \eqref{Mod:CovOp}. This implies that the random field $Z_{k} (\mathbf{s})$ has a nugget effect $ \langle C_{\varepsilon}(\mathbf{s},\mathbf{s}; \eta_{k}), \eta_{k} \rangle$ that should be considered when fitting a covariance (variogram) model. 

The covariance estimator of $Y(\mathbf{s};v)$ is obtained from \eqref{kernelEst}, after estimating the covariance $\boldsymbol{\Sigma}(h)= \mathrm{Cov} \{\mathbf{Z}(\mathbf{s}_{1}), \mathbf{Z}(\mathbf{s}_{2}) \}$, $h= \|\mathbf{s}_{1} - \mathbf{s}_{2}  \|$, from the data $\{ \langle y(\mathbf{s}_{i};\cdot), \eta_{k} \rangle \}_{i=1}^{n}$, $k=1, \ldots, K$. These ideas are extended to estimating the variogram of $Y(\mathbf{s};v)$ by estimating the matrix-variogram of $\mathbf{Z}(\mathbf{s})$. 

In the case of the trace-variogram of $Y(\mathbf{s};v)$,  from \eqref{small-v-m} we have that the mean of the inner product satisfies 
 $
\mathbb{E} \{ \langle Y(\mathbf{s}_{i};\cdot),  Y(\mathbf{s}_{j};\cdot) \rangle  \} = \mathbb{E} \{ \langle X(\mathbf{s}_{i};\cdot),  X(\mathbf{s}_{j};\cdot ) \rangle  \} + \mathbb{E} \{ \langle \varepsilon(\mathbf{s}_{i};\cdot),  \varepsilon(\mathbf{s}_{j};\cdot ) \rangle  \}
$. That is, the trace-covariogram of $Y(\mathbf{s};v)$ is such that
$$
 \sigma_{\mathrm{tr,Y}}(\mathbf{s}_{1}, \mathbf{s}_{2} )=  \sigma_{\mathrm{tr,X}}(\mathbf{s}_{1}, \mathbf{s}_{2} ) + \mathds{1} ( \mathbf{s}_{1}=\mathbf{s}_{2}) \sigma_{\mathrm{tr, \varepsilon}}(\mathbf{s}_{1}, \mathbf{s}_{2} ).
$$
Then, when estimating  $\sigma_{\mathrm{tr,Y}}(\mathbf{s}_{1}, \mathbf{s}_{2} )$ using $\{ \langle y(\mathbf{s}_{i};\cdot), y(\mathbf{s}_{j};\cdot) \rangle \}_{i,j}$, as described in Section~\ref{EstCov}, one should consider the nugget effect  $\sigma_{\mathrm{tr, \varepsilon}}(\mathbf{s}, \mathbf{s} )$. The same is true in the case of the trace-variogram.

\subsection{Large-scale and Small-scale Variations}\label{section:mix-v}

Datasets often have both a trend component (large-scale variation) and a spatial variability (small-scale variation), e.g., temperature data show  an increasing tendency and a small variability around this tendency. In the context of functional data correlated only in time,  \cite{MtzHdz2019} proposed a method to estimate trend using tensor product surfaces. 

A model for spatial functional data can be written as 
\begin{equation*}
Y(\mathbf{s};v) =\mu(\mathbf{s};v) + X(\mathbf{s}; v) + \varepsilon(\mathbf{s};v).
\end{equation*}
As before, the mean can be expressed as $\mu(\mathbf{s};v)= \sum_{l=0}^{L} a_{l}( \mathbf{s}) f_{l}(v)$. Similarly to \eqref{LD-model}, the model to estimate the parameters of the mean can be written as
\begin{equation}\label{meanModel-LS}
Y(\mathbf{s};v) = \sum_{l=0}^{L} a_{l}( \mathbf{s}) f_{l}(v) + \epsilon(\mathbf{s};v),
\end{equation}
where the residual $\epsilon(\mathbf{s};v) :=  X(\mathbf{s}; v) + \varepsilon(\mathbf{s};v)$ is now a functional random field with mean zero and covariance $C_{\epsilon} (\mathbf{s}_{1}, \mathbf{s}_{2}; \sqcdot)$. Thus, unlike model  \eqref{LD-model}, in which the residuals are not correlated, model \eqref{meanModel-LS} has spatially correlated residuals. Because of this correlation, we can use the GLS method instead of the OLS method \citep[][]{menafoglio2013}. Let $\boldsymbol{\Sigma}$ be the trace-variogram matrix of $\epsilon(\mathbf{s}_{i};v)$ at different distances of the locations $\mathbf{s}_{i}$. Then, the GLS estimator of $\mu(\mathbf{s}_{i};v)$ is
\begin{equation}\label{GLD-mean}
\hat{\boldsymbol{\mu}}(v)= \mathbf{A} (\mathbf{A}^{\mathsf{T}}\boldsymbol{\Sigma}^{-1}\mathbf{A})^{-1} \mathbf{A}^{\mathsf{T}}\boldsymbol{\Sigma}^{-1} \mathbf{y}(v),
\end{equation}
where $\mathbf{A}$ is the design matrix, and $ \mathbf{y}(v)$ the evaluation of $y(\mathbf{s}; v)$ at $n$ locations, both defined in Section \ref{section:large-v}. 

The  trace-variogram matrix $\boldsymbol{\Sigma}$ cannot be estimated directly from $\epsilon(\mathbf{s}_{i};v)$, because we do not observe them. The approach commonly used, to solve this problem, is to use the following  iterative procedure: 
\begin{enumerate}
\item Compute the initial estimator $\hat{\mu}_{0}(\mathbf{s}_{i};v)$ as in  \eqref{estimatorLD}, assuming that $\{\epsilon(\mathbf{s}_{i};v)\}_{i=1}^{n}$ are not spatially correlated. 
\item Compute the residuals $\hat{\epsilon}(\mathbf{s}_{i};v)= y(\mathbf{s}_{i};v)- \hat{\mu}_{0}(\mathbf{s}_{i};v)$, for $i=1,\ldots,n$.
\item Estimate the initial empirical trace-variogram $\hat{\gamma}_{\mathrm{tr},0}(h)$ as in \eqref{EmTraVar}, using $\{\hat{\epsilon}(\mathbf{s}_{i};v)\}_{i=1}^{n}$, and then, obtain an estimator $\hat{\boldsymbol{\Sigma}}_{0}$ of $\boldsymbol{\Sigma}$ by fitting a parametric model with a nugget effect.
\item Re-estimate  the mean to get $\hat{\mu}_{1}(\mathbf{s}_{i};v)$ by using \eqref{GLD-mean} with $\hat{\boldsymbol{\Sigma}}_{0}$.
\item Repeat steps $2-4$ until convergence.
\end{enumerate}
Once the mean $\mu(\mathbf{s};v)$ is estimated, it is removed from the data. Then, the covariance is estimated as in Section \ref{section:small-v}, using the spatial functional data $\{ y(\mathbf{s}_{i};v) - \hat{\mu}(\mathbf{s}_{i};v)\}_{i=1}^{n}$.   

The iterative procedure can also be performed with the finite basis function $\boldsymbol{\eta}(v)$. In step $1$, we can use \eqref{estimatorLD2}, and in step $4$, use the estimator $\mathrm{vec}(\hat{\mathbf{B}})= ( \mathbf{J}_{\boldsymbol{\eta}}^{\mathsf{T}}\otimes \mathbf{A}^{\mathsf{T}}\hat{\boldsymbol{\Sigma}}_{0} \mathbf{A} )^{-1}  \mathrm{vec}(\mathbf{A}^{\mathsf{T}}\hat{\boldsymbol{\Sigma}}_{0} \mathbf{Z} \mathbf{J}_{\boldsymbol{\eta}})$  to obtain $\hat{\mu}_{1}$.

\section{Kriging for Functional Random Fields}\label{KrigingSFD}

In spatial statistics, the concept of kriging (co-kriging for the multivariate setting) is a synonym of optimal interpolation. The main goal is to be able to predict at locations where data are not observed. This predictor is a linear combination of the observed data, such that it is the best linear unbiased predictor under squared loss. Here, we briefly mention the concept of kriging and we redirect readers to the references provided. 

Let $\mathbf{s}_{0}\in D$ be the location at which the curve will be predicted. Let $\{\Psi_{1}, \ldots, \Psi_{n} \}$ be linear operators from $\mathcal{H}$ to $\mathcal{H}$. In general, kriging can be defined as 
\begin{equation}\label{kriging}
\hat{x}(\mathbf{s}_{0}; v )= \sum_{i=1}^{n} \Psi_{i} \{x(\mathbf{s}_{i}; v) \},
\end{equation}
where the coefficients $\Psi_{i}$ are obtained by minimizing the square norm of the error prediction, $\hat{x}(\mathbf{s}_{0}; v )- X(\mathbf{s}_{0}; v )$, with an additional constraint of unbiasedness. That is,   
\begin{align*}
\min_{\Psi_{1},\ldots, \Psi_{n} } & \mathbb{E}\{ \langle \hat{x}(\mathbf{s}_{0}; \cdot )- X(\mathbf{s}_{0}; \cdot ), \hat{x}(\mathbf{s}_{0}; \cdot )- X(\mathbf{s}_{0}; \cdot ) \rangle \} \\
\mbox{ s.t. } & \mathbb{E} \{ \hat{x}(\mathbf{s}_{0}; v )\}=\mathbb{E} \{ X(\mathbf{s}_{0}; v )\} .
\end{align*}
A particular case of the coefficients operators $\Psi_{i}$ in \eqref{kriging} is the so-called \textit{kernel operators}, which are defined as $\Psi_{i} (f)(v)= \int_{0}^1 \! \lambda_{i} (v,u) f(u) \mathrm{d}u $, $f\in \mathcal{H}$. In this case, the estimation is through functions $\lambda_{i} (v,u)$. Other cases are $\Psi_{i} (f)(v)= \lambda_{i} (v) f(v)$ and $\Psi_{i} (f)(v)= w_{i} f(v)$, with $w_{i}$ scalars. The latter corresponds to a simple ponderation of the observed curves, i.e., $\hat{x}(\mathbf{s}_{0}; v )= \sum_{i=1}^{n} w_{i} x(\mathbf{s}_{i}; v)$. All these cases can be fitted  into ordinary or universal kriging. For ordinary kriging, see \cite{Goulard1993},  \cite{NERINI2010}, \cite{Giraldo2010}, and \cite{GiraldoEtAl2011}, as well as the review paper by \cite{DelicadoEtAl2010}. For universal kriging, see \cite{menafoglio2013}, \cite{Caballero2013}, \cite{ReyesEtAl2015}, and \cite{MENAFOGLIO201639}.  For co-kriging (multivariate functional random fields), see  \cite{Bohorquez2017} and \cite{Grujic2018}. An alternative approach to kriging, which is based on a tensor function space,  can be found in \cite{AguileraMorillo2017}.

\section{Surface Time Series}\label{SurfaceD}

In many phenomena, data can be collected in the form of a surface, called \textit{surface data} in our case, or \textit{manifold data} for more complex structures. For example, we can have data arising from neuroimaging \citep[][]{LilaEtAl2016}, from two-dimensional time-frequency domains \citep[][]{AstonEtAl2017}, from satellite images  \citep[][]{Zhang2011}, and  functional data with two-dimensional domain \citep[][]{CrainiceanuEtAl2011,MorrisEtAl2011}. Here, we consider surface data as functional,  that is, the atoms of the functional random variable are continuous surfaces. Surface data provide an alternative approach to analyzing spatial data, where the continuous  realization of a random field is considered as a unit. This approach can have computational advantages, especially if the locations, where data are observed, are dense in space. Particularly, spatio-temporal data \citep[][]{CressieEtAl2011} can be considered as surface data that are observed over time (surface time series). Surface data capture the  spatial dependence through the continuity of the surface \cite[see, e.g.,][]{Bernardi2017}. Moreover, the approach of surface data  can be applied to a nonplanar spatial domain, such as a sphere or  a general two-dimentional manifold   \citep[][]{DassiEtAl2015,Wilhelm2016,EttingerEtAl2016,MENAFOGLIO2017401,GrecoEtAl2018}. Kriging can be applied to these complex domains by using an appropriate  distance, but the covariance models do not necessarily guarantee a positive definite covariance, e.g., the Mat\'ern covariance \citep[][]{Gneiting2013}. Here, we focus on spatio-temporal data, where at each time point we observe a surface. 

The functional time series approach to spatial statistics has been studied by  \cite{Ruiz-MedinaEtAl2007}. In \cite{AstonEtAl2017}, a tensor product Hilbert space was considered to propose a separability test for the covariance operators of random surfaces.

\subsection{Basic Concepts}

Let $(\Omega, F, P)$ be a probability space, and let  $\mathcal{H}$ be the  Hilbert space defined as the set of functions with  domain $D\subset \mathbb{R}^{2}$, $\mathcal{H}= \{ f:D\to \mathbb{R}: \int_{D}  | f(\mathbf{s} )|^{2} \mathrm{d}\mathbf{s} <\infty \}$, and with the inner product $\langle f,g \rangle=\int_{D}  f(\mathbf{s}) g(\mathbf{s} ) \mathrm{d}\mathbf{s}$. The norm induced by the inner product is denoted by $\| \cdot \|_{\mathcal{H}}$. Thus, a random variable $X: \Omega \to \mathcal{H}$ is a functional random variable with a surface as atom. We denote by $X(\mathbf{s})$ this functional random variable with $\mathbf{s}= (s_{1}, s_{2})\in D$.

A functional time series is a sequence of functional random variables $\{X_{t}(\mathbf{s}); t\in \mathbb{Z}\}$ in $\mathcal{H}$. \cite{Bosq2000} 
is a monograph on linear processes in function spaces, including functional time series in a Hilbert space.  \cite{HORMANN2012157} reviewed  functional time series. Also, see \cite{Ramsay-Silverman2005}, \cite{FerratyVieu2006}, and \cite{HorvathKokoszka2012} for a general introduction to functional data analysis.
In this paper,  we assume that $ \mathbb{E} \{ \| X_{t}(\mathbf{s}) \|_{\mathcal{H}}^{2} \} < \infty$ for all $t\in \mathbb{Z}$.

The mean of the surface time series $\{X_{t}(\mathbf{s})\}$ is defined as $\mu_{t}(\mathbf{s})= \mathbb{E}\{X_{t}(\mathbf{s})\}$, where $\mu_{t}(\mathbf{s}) $ is such that $\mathbb{E}(\langle X_{t}, f \rangle )= \langle \mu_{t}, f \rangle $ for all $f\in \mathcal{H}$. 
The covariance function at lag $h\in \mathbb{N}$ is defined as  
$$C_{X_{t-h}, X_{t}}(f)= \mathbb{E} \{ \langle X_{t-h}- \mu_{t-h}, f \rangle  ( X_{t}-\mu_{t} ) \}, \quad f\in \mathcal{H}.$$
This covariance function can be expressed as 
$$C_{X_{t-h}, X_{t}}(f)(\cdot)= \int \sigma_{t-h, t}( \cdot,\mathbf{s} ) f(\mathbf{s}) \mathrm{d} \mathbf{s}, $$
where $\sigma_{t-h, t}( \mathbf{s}_{1},\mathbf{s}_{2} )= \mathrm{Cov} \{  X_{t-h}( \mathbf{s}_{1}),  X_{t}(\mathbf{s}_{2}) \}$.  
The stationarity condition is important for statistical inference.
\begin{Def}[Weak stationarity] A surface time series $\{X_{t} (\mathbf{s}); t\in \mathbb{Z}\}$ is said to be (weakly) stationary if 
\begin{enumerate}
\item $\mathbb{E} \{ X_{t} (\mathbf{s}) \}= \mu(\mathbf{s}) $ for all $t\in \mathbb{Z}$, and
\item $C_{X_{t_{1}+h}, X_{t_{2}+h}}(f) = C_{X_{t_{1}}, X_{t_{2}}}(f)$ for all $t_{1}, t_{2}  \in \mathbb{Z}$, $h\in \mathbb{N}$, and $f\in \mathcal{H}$.
\end{enumerate}
\end{Def}	
If the surface time series is stationary, we write $C_{h}$ for the covariance functions instead of $C_{X_{t}, X_{t+h}}$. The definition of stationarity does not require stationarity over the space $D$, e.g., for each $t$, data can be a realization of a nonstationary random field. In general, the covariance function $C_{0}$ describes only the spatial dependence, whereas $C_{h}$, for $h\neq 0$, describes the dependency over time of the surface time series.

A \textit{surface white noise} is a stationary surface time series with a zero mean and a covariance function $C_{h}=0$, if $h\neq 0$. Thus,  the surface white noise can have a spatial correlation at each time point, but not across time.

Similarly as before, the eigenfunctions are defined as functions $ \zeta \in \mathcal{H}$ such that 
\begin{equation*}
C_{0} (\zeta)  (\mathbf{s}) = \lambda \zeta (\mathbf{s}),
\end{equation*}
where $\lambda$ is positive and is the corresponding eigenvalue.  Moreover, the covariance operator $C_{0}$ can be decomposed in terms of the eigenfunctions, that is 
$$
C_{0} (f)(\mathbf{s}) =  \sum_{j=1}^{\infty} \lambda_{j} \langle   \zeta_{j}, f \rangle \zeta_{j} (\mathbf{s}),
$$
where $\zeta_{j}$, $j=1,2,\ldots$, are the eigenfunctions of $C_{0}$ with eigenvalues $ \lambda_{j}$. The eigenvalues are such that $\sum_{j=1}^{\infty} \lambda_{j} = \mathbb{E} \{\|X_{0}(\mathbf{s}) \|^{2}_{\mathcal{H}} \}< \infty$. This operator $C_{0}$ is nuclear and therefore Hilbert-Schmidt. Similarly as in the finite dimensional case, eigenfunctions are important to reduce dimensionality.

\subsection{Estimation}

Now, we describe estimators of the mean $\mu(\mathbf{s})$ and the covariance operator $C_{h}$. Let $\{ x_{t}( \mathbf{s}) \}_{t=1}^{T}$ be a realization of a stationary surface time series $X_{t}( \mathbf{s})$ with mean $\mu(\mathbf{s})$. The sample mean is defined as 
\begin{equation}\label{SD-mean}
\hat{\mu}( \mathbf{s})= \frac{1}{T} \sum_{t=1}^{T} x_{t}( \mathbf{s}).
\end{equation}
The sample mean $\hat{\mu}( \mathbf{s})$ is an unbiased estimator of $\mu( \mathbf{s})$.

The empirical covariance at lag $h$ of $X_{t}( \mathbf{s})$ is defined as 
\begin{equation}\label{Cov-FAR-Est}
\hat{C}_{h}(f)(\mathbf{s}) = \frac{1}{T-h} \sum_{t=1}^{T-h} \langle  x_{t}- \hat{\mu} , f \rangle \{ x_{t+h}(\mathbf{s}) - \hat{\mu}(\mathbf{s}) \}, \quad f\in \mathcal{H},
\end{equation}
and the corresponding empirical kernel is defined as 
$$
\hat{\sigma}_{h} (\mathbf{s}_{1}, \mathbf{s}_{2} ) =  \frac{1}{T-h} \sum_{t=1}^{T-h}  \{ x_{t+h}(\mathbf{s}_{1}) - \hat{\mu}(\mathbf{s}_{1}) \}  \{ x_{t}(\mathbf{s}_{2}) - \hat{\mu}(\mathbf{s}_{2}) \}.
$$
The empirical covariance operator $\hat{C}_{h}$ is an unbiased estimator of  $C_{h}$, see \cite{Bosq2000}. For papers related to the mean and the covariance functions, refer to previous studies by authors \cite{HormanKokoszka}, \cite{HorvathEtAl2013}, and \cite{HorvathEtAl2014}.

In time series analysis, there is no direct modeling of the covariance function. Instead, a model for the process is proposed, and the covariance function is derived from the model. We adopt this idea in the next section (Section \ref{Sec:STSModeling}). 

\subsection{Modeling}\label{Sec:STSModeling}
This section discusses two topics: the continuous estimation of the surface and the modeling of the continuous surface series.

\subsubsection{Estimating the continuous surface}\label{Est-CS}
In practice, data are observed on a finite set of points, i.e., for each $t=1, \ldots, T$, we observe $n_{t}$ points $\{y_{t,i}= y_{t}(\mathbf{s}_{i})\}_{i=1}^{n_{t}}$ of the functional data $Y_{t}(\mathbf{s})$ on a set of points $\{\mathbf{s}_{1}, \ldots, \mathbf{s}_{n_{t}}\}\subset D$, and possibly with measurement errors. Thus, for each $t$, the unknown surface (deterministic field) $y_{t}(\mathbf{s})$ needs to be estimated. This continuous surface estimate can be associated with kriging in classical spatial data analysis.

Because the procedure to estimate  $y_{t}(\mathbf{s})$  is independent of $t$, we drop the subindex $t$ in the sequel, and we consider $n$ as the sample size. Thus, a model of $y(\mathbf{s})$ can be written as 
\begin{equation*}
y_{i}= y(\mathbf{s}_{i})+ \varepsilon_{i},
\end{equation*}
where $\{\varepsilon_{i}\}_{i=1}^{n}$ represent the measurement errors that are spatially uncorrelated. The function $y(\mathbf{s})$ describes the spatial structure of the phenomenon being studied.

To estimate $y(\mathbf{s})$, one can extend the smoothness techniques of the curves described in \cite{Ramsay-Silverman2005}  to surfaces (manifolds). In particular, one can extend the spline smoothing.  The extension of spline smoothing to surfaces is an important research area. Extensions have been done over Euclidean domains and non-Euclidean domains, including the spherical domain. One extension is to use the  tensor product of univariate B-splines \citep[][]{EilersMarx96,Wood2006,QingguoEtAl2010,XiaoEtAl2013}.  In this case, an estimator of $y (\mathbf{s})$ has the form 
 $$\hat{y}( \mathbf{s}) = \sum_{k=1}^{K_{1}} \sum_{l=1}^{K_{2}} \theta_{kl} \eta_{k}(s_{1}) \nu_{l}(s_{2}),$$
 where $\{\eta_{k}\}_{k=1}^{K_{1}}$ and $\{\nu_{l}\}_{l=1}^{K_{2}}$ are B-splines basis functions for $s_{1}$ and $s_{2}$ coordinates, respectively, with $\mathbf{s}=(s_{1}, s_{2})$. In the estimation procedure, smoothness properties are imposed through a penalization term. To read about some approaches of bivariate smoothing, see  \cite{Ruppert2003}, \cite{lai2007}, and \cite{Wood-Book}.

In general, the estimation of $y(\mathbf{s})$ can be formulated as the minimization of the sum of squared errors with a penalization term. The penalization term measures the roughness of the fitted surface and can carry partial information of $y(\mathbf{s})$. That is,  the optimization  problem is written as  
\begin{equation}\label{op:sp}
\sum_{i=1}^{n} \{ y_{i}- y(\mathbf{s}_{i}) \}^{2} + \lambda P(y),
\end{equation}
where $P(y)$ is the penalization term, and $\lambda$ the smoothness parameter which controls the smoothness of the estimated surface. A popular penalization is the thin-plate energy, which is defined as $P(y)= \int \{( \frac{\partial^{2} y}{\partial s^{2}_{1}})^{2} + 2 (\frac{\partial^{2} y}{\partial s_{1} \partial s_{2}})^{2} + (\frac{\partial^{2} y}{\partial s^{2}_{2}} )^{2} \} \mathrm{d} \mathbf{s}$. The resulting estimator is called the thin plate splines \citep[][]{Duchon77}. A Bayesian adaptive thin plate spline was proposed in \cite{YueAtAl2010}.  Another example of  penalization involves the Laplacian, that is $P(y)= \int (  \frac{\partial^{2} y}{\partial s^{2}_{1}} +  \frac{\partial^{2} y}{\partial s^{2}_{2}} )^{2} \mathrm{d} \mathbf{s}$ \citep[][]{WoodEtAl2008,SangalliEtAl2013}. The definition of the penalty term depends on each specific problem. For example, in the Laplacian case, the unique penalty parameter $\lambda$ that  controls both directions $s_{1}$ and $s_{2}$ implies an isotropic smoothing. In contrast, if $P(y)= \lambda_{1} P_{1} (y) + \lambda_{2}P_{2} (y)$, where $P_{i} (y)$ is a penalty term in the $i$-th coordinate, then it results in anisotropic smoothing.

In general, the penalty term can be defined in terms of a partial differential equation (PDE). For example, $P(y)= \int (Ly -u )^{2} \mathrm{d} \mathbf{s}$, where $L$ is a differential operator and $Ly=u$ is a PDE \citep[see, e.g.,][]{AzzimonteEtAl2015,Sangalli2019}. The advantage of the penalty term with PDE  is that it can  handle complex domains with boundary conditions or  interior holes, and varies  depending on the phenomena being  studied.  The PDE is such that it contains information about the phenomena, and it regularizes the estimation with values of $\lambda$.  The solution of \eqref{op:sp} may not have a closed form, but it can be approximated by using finite elements  analysis, see \cite{Ramsay2002}, \cite{DuchampEtAl2003}, and \cite{SangalliEtAl2013}.

We describe the solution of \eqref{op:sp} using the finite elements analysis technique. Let $\mathcal{M}$ be a mesh of $D$. Let $\{\phi_{1}(\mathbf{s}),\ldots,\phi_{K}(\mathbf{s})\}$ be basis functions that are piece-wise polynomials associated with the mesh $\mathcal{M}$. Then, the estimator of $y(\mathbf{s})$ is assumed to have the form 
\begin{equation*}
\hat{y}(\mathbf{s}) = \sum_{k=1}^{K} \beta_{k} \phi_{k}(\mathbf{s}), 
\end{equation*}
where the coefficients $\boldsymbol{\beta}= (\beta_{1}, \ldots, \beta_{K} )^{\mathsf{T}}$ need to be estimated. Let $\mathbf{y}= ( y_{1}, \ldots, y_{n} )^{\mathsf{T}}$ be the observed values over $D$, and let $\mathbf{P}$  be the discretization of the penalty. Then, the estimator of $\boldsymbol{\beta}$ has the form
\begin{equation*}
\hat{\boldsymbol{\beta}}= \left(\varPhi^{\mathsf{T}}  \varPhi + \lambda \mathbf{P} \right)^{-1} \varPhi^{\mathsf{T}}  \mathbf{y},
\end{equation*}
where $\varPhi = \{\phi_{k} (\mathbf{s}_{i}) \}_{i,k=1}^{n,K}$ is the $n\times K$ matrix which represents the evaluation  of each basis function at the locations at which data are observed.

The time component can also be considered in the PDE. \cite{ARNONE2019275} proposed general forms of time-dependent PDEs in the context of  time dependent surface data.

Additional studies of data over complex domains have been published. \cite{Wang2007} proposed a modified thin plate spline over a complex domain; \cite{LindgrenAndHavard2011} linked Gaussian fields via stochastic partial differential equations, where the solution is found using a finite elements analysis; \cite{Scott-Hayward2014} proposed a complex region spatial smoother using the geodesic distance; and  \cite{Menafoglio2018} proposed a methodology for spatial fields of object data over complex domains.

\subsubsection{Functional autoregressive models}

Here, we assume that the functional time series consist of continuous surfaces (deterministic fields) that can be estimated with the methods described in Section \ref{Est-CS}. In the context of functional time series, the most popular model is the functional autoregressive model of order $P$, FAR$(P)$. A surface time series $\{X_{t}(\mathbf{s})\}$ follows the FAR$(P)$ model if $X_{t}(\mathbf{s}) = \sum_{p=1}^{P} \Psi_{p} ( X_{t-p})(\mathbf{s}) + W_{t}(\mathbf{s})$, where each coefficient $\Psi_{p}: \mathcal{H} \to \mathcal{H}$ is an operator, and $\{W_{t}(\mathbf{s})\}$ is a surface white noise. In practice, the order $P$ needs to be estimated \citep[][]{KokoszkaR2013}. Here, we assume $P=1$ to illustrate the ideas  and to simplify the notations.

 Let $\{ y_{t}(\mathbf{s}) \}_{t=1}^{T}$ be the surface data observed over time $t=1, \ldots, T$, and assume that it is a stationary surface time series. Then, the dependency over time can be modeled by using the FAR$(1)$ process, that is
\begin{align}
y_{t}(\mathbf{s}) & = \mu(\mathbf{s}) + x_{t}(\mathbf{s}), \label{SFTS1} \\
x_{t}(\mathbf{s}) &  = \Psi (x_{t-1}  )(\mathbf{s}) + W_{t} (\mathbf{s}),  \label{SFTS2} 
\end{align}
where $\mu(\mathbf{s})$ represents the surface mean, i.e., the large-scale variation for all $t$. The unobserved $x_{t}(\mathbf{s})$ follows a stationary FAR$(1)$ process with mean zero, and $\{ W_{t} (\mathbf{s})\}$ is a surface white noise. The surface white noise $\{W_{t} (\mathbf{s}) \}$ can be interpreted as the surface data components that describe the small-scale variation for each time $t$,  which are not correlated over time. The dependency over time is driven by the operator $\Psi$. 

Since $\{y_{t}(\mathbf{s})\}$ is assumed to be stationary, the estimation of $\mu(\mathbf{s})$ can be obtained as in \eqref{SD-mean}. After removing the mean, the rest of the analysis is performed on the process $x_{t}(\mathbf{s}) = y_{t}(\mathbf{s}) - \hat{\mu}_{t}(\mathbf{s})$.

Now, we focus on the estimation of the coefficient operator $\Psi$. Let $C_{h}$ be the covariance operator of the FAR$(1)$ process $X_{t} (\mathbf{s})$ with realization $\{x_{t} (\mathbf{s})\}$.  Then, the  covariance operator of $X_{t}(\mathbf{s})$ satisfies 
\begin{equation}\label{Cov0-FAR}
C_{1}(f)= \Psi \{ C_{0} (f)\}, \quad f\in \mathcal{H}.
\end{equation}
Moreover, for $h\in \mathbb{N}$, we have that $ C_{h}(f)= \Psi^{h} \{ C_{0} (f)\}$. So if the coefficient $\Psi$ and $C_{0}$ are known, then we can compute the covariance function at any lag $h\in \mathbb{N}$. Thus, the estimation of $\Psi$ is crucial.

To obtain an estimator of $\Psi$, we can use the estimation of $C_{0}$ and $C_{1}$ in  \eqref{Cov0-FAR}, defining  $\Psi (f) = C_{1} \{ C_{0}^{-1} (f)\}$, if $C_{0}$ is invertible. In principle, we can always use the estimator \eqref{Cov-FAR-Est} for $h=0,1$, and compute the inverse $\hat{C}^{-1}_{0}$. However, when the sample size $T$ tends to infinity, $\hat{C}^{-1}_{0}$ becomes unbounded \citep[][]{CardotFerratySarda99}. This is because $C_{0}$ is a compact operator \citep[][]{Bosq2000}. Thus, it is necessary to use some regularization methods to obtain $C_{0}^{-1}$. That is, $(C_{0} + \alpha_{T})^{-1}$ is computed instead of $C_{0}^{-1}$ where $\alpha_{T}>0$ and $\alpha_{T} \downarrow 0$.  Alternatively, $C_{0}^{-1}$ can be approximated by using only the first $k$ eigenfunctions corresponding to the largest eigenvalues, that is, $C_{0}^{-1}= \sum_{j=1}^{k} \lambda_{j}^{-1} \zeta_{j} \otimes \zeta_{j} $, see \cite{Bosq2000} and \cite{KokoszkaEtAl2017}.  Let $\hat{C}_{0}^{-1}$ denote an estimator of the inverse operator, either using some regularization method, finite eigenfunctions or other methods \citep[][]{MARTINEZHERNANDEZ201966}. Then, the estimator of $\Psi$ is defined as 
$$
\hat{\Psi}(f)= \hat{C}_{1} \{\hat{C}_{0}^{-1} (f) \}, \quad f\in \mathcal{H}.
$$
Once $\mu(\mathbf{s})$ and $\Psi$ are estimated, then the one-step-ahead prediction is obtained as 
$$
\hat{y}_{T+1}(\mathbf{s}) =\hat{\mu}( \mathbf{s})+  \hat{x}_{T+1}  (\mathbf{s}),
$$
where $\hat{x}_{T+1}(\mathbf{s}) = \hat{\Psi} (x_{T} ) (\mathbf{s}).$ The estimators of $\mu$, $\Psi$, and $C_{h}$ can be explicitly expressed in terms of the basis functions $\{\phi_{1}(\mathbf{s}),\ldots, \phi_{K}(\mathbf{s}) \}$ from the finite elements technique.  Alternative approaches can be used to predict data, for example one can extend the ideas described in \cite{Hyndmanetal2007} and \cite{AueEtAl2015} to surfaces. 

The approach of the surface time series provides alternatives to analyze spatio-temporal data. The advantage of this approach is that it can handle  data  collected on a large scale for each time point, over a general domain, Euclidean, or non-Euclidean. With this approach, the process is modeled instead of the covariance, which is convenient when the classical covariance models are not guaranteed to be  positive definite. 

\section{Statistical Software}\label{Software}

Here, we mention some packages available in \texttt{R} for SFD and surface data.  The package \textit{fda} \citep[][]{fda} provides several commands to analyze and construct continuous functions. It contains several options of basis functions, such as Fourier basis functions and  spline basis functions. The basis functions can be used to estimate the continuous functions as in \eqref{approxX}. 
Once the coefficients are obtained in \eqref{approxX}, we can use packages for classical spatial data, e.g., \textit{spatial} \citep[][]{VenablesEtAl2002}, \textit{gstat} \citep[][]{gstat}, \textit{RandomFields} \citep[][]{RandomField}, \textit{fields} \citep[][]{fields}, \textit{geoR}  \citep[][]{geoR}, \textit{ExaGeoStatR} \citep[]{Abdulah2019} for large datasets, and \textit{spBayes} \cite[][]{spBayes} for Bayesian analysis of hierarchical multivariate models. Thus, when data are expressed in terms of basis functions as in \eqref{approxX}, we can combine the \textit{fda} package and the packages for spatial data  to obtain estimators of the mean and the covariance functions described in Section \ref{sec:Est}.

The package \textit{geofd} by \cite{geofd} implements kriging of functional data described in Section \ref{KrigingSFD}. The curves observed are pre-processed by fitting Fourier or B-splines basis functions. Also, this package provides a command to compute the trace-variogram  defined in \eqref{EmTraVar}. Another package related to Section \ref{KrigingSFD} is \textit{fdagstat} \citep[][]{fdagstat}. This package  implements kriging, cokriging, and universal kriging, and includes the large-scale variation described in Sections \ref{section:large-v} and  \ref{section:mix-v}.

For the surface data described in Section \ref{SurfaceD}, one can use the package \textit{mgcv} by \cite{Wood-Book}, which allows us to smooth surfaces. The package \textit{fdaPDE} \citep[][]{fdaPDE}  implements smoothing with PDE penalization described in \cite{SangalliEtAl2013} and \cite{AzzimonteEtAl2015}. \textit{INLA} \cite[][]{RueEtAl2009} can be used to estimate continuous surfaces.
 The package \textit{Manifoldgstat} by \cite{Manifoldgstat} implements kriging for manifold-valued random fields. 

Some visualization tools for functional data and functional time series are the functional boxplots \citep[][]{SunGenton2011,SunGenton2012} implemented in the \textit{fda}  package; for functional images and surfaces, surface boxplots are used \citep[][]{GentonEtAl2014}. These tools are based on an order induced by a depth notion for functional data.

\section{Data Analysis}\label{DataA}
In this section, we illustrate general ideas of  modeling spatial functional data and surface data, without a deep statistical analysis of the data. Our goal is to provide a general example, using only available packages. 

We use wind data simulated from the Weather Research and Forecasting (WRF) model by \cite{Yip2018}. Each measurement corresponds to hourly wind speed from $2009$ to $2014$, in a $115 \times 115 $ km region centered in Dumat Al Jandal, Saudi Arabia. That is where the first wind farm of the country is being built.  These wind speeds  are simulated on a regular grid of points in space, namely at $5$-km resolution.

\subsection{Spatial Functional Data Approach} 

Our data are the daily wind speed, where  $y(\mathbf{s}_{i};v_{j})$ is the wind speed at location $\mathbf{s}_{i}$ for hour $v_{j}$, for  $v_{j}= 1,2,\ldots, 24$. With the functional approach, we consider $y(\mathbf{s}_{i};\cdot)$ as a single object, assuming continuity over time. 
\begin{figure}[!t]
\begin{center}
\includegraphics[scale=.47]{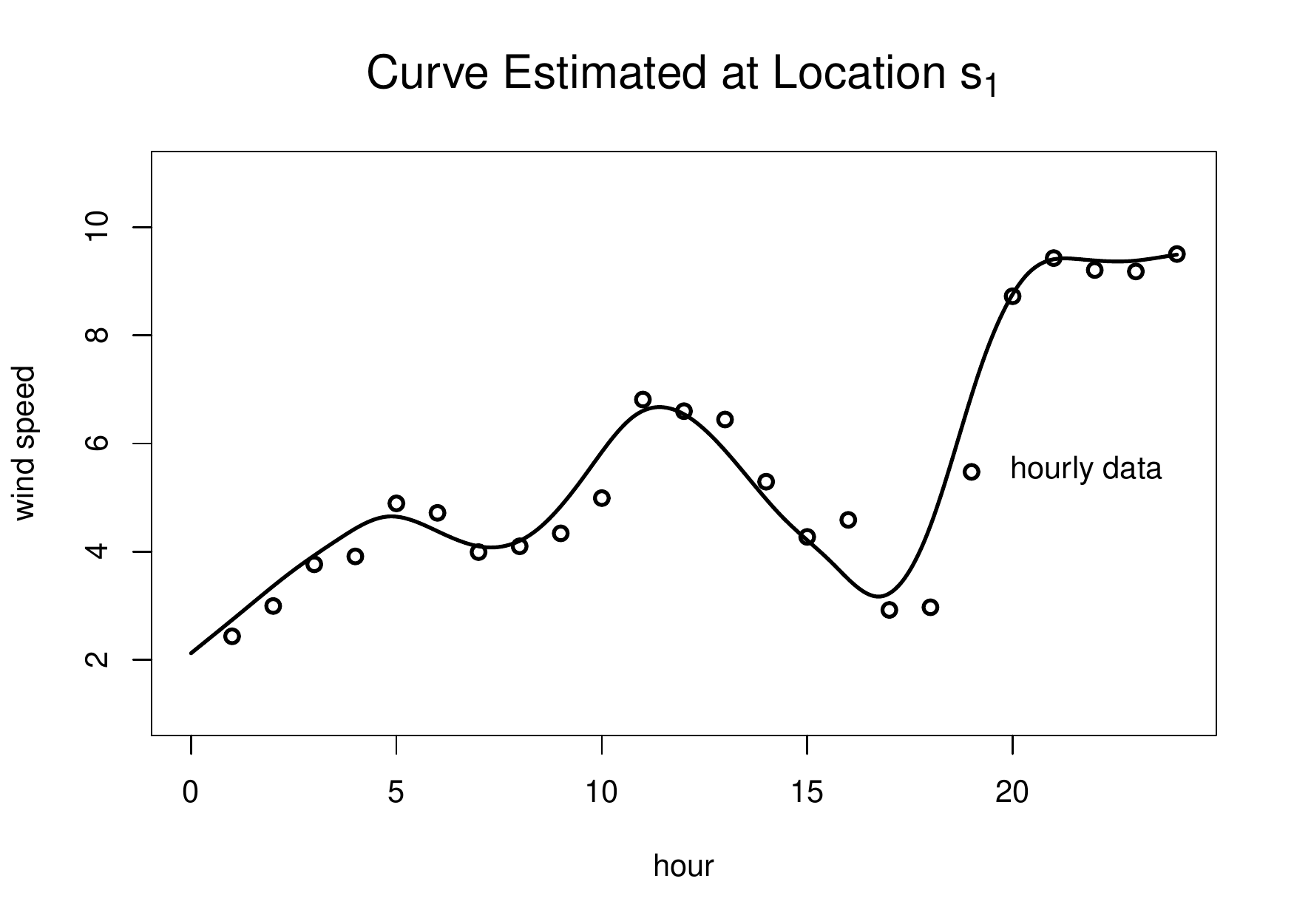} 
\includegraphics[scale=.47]{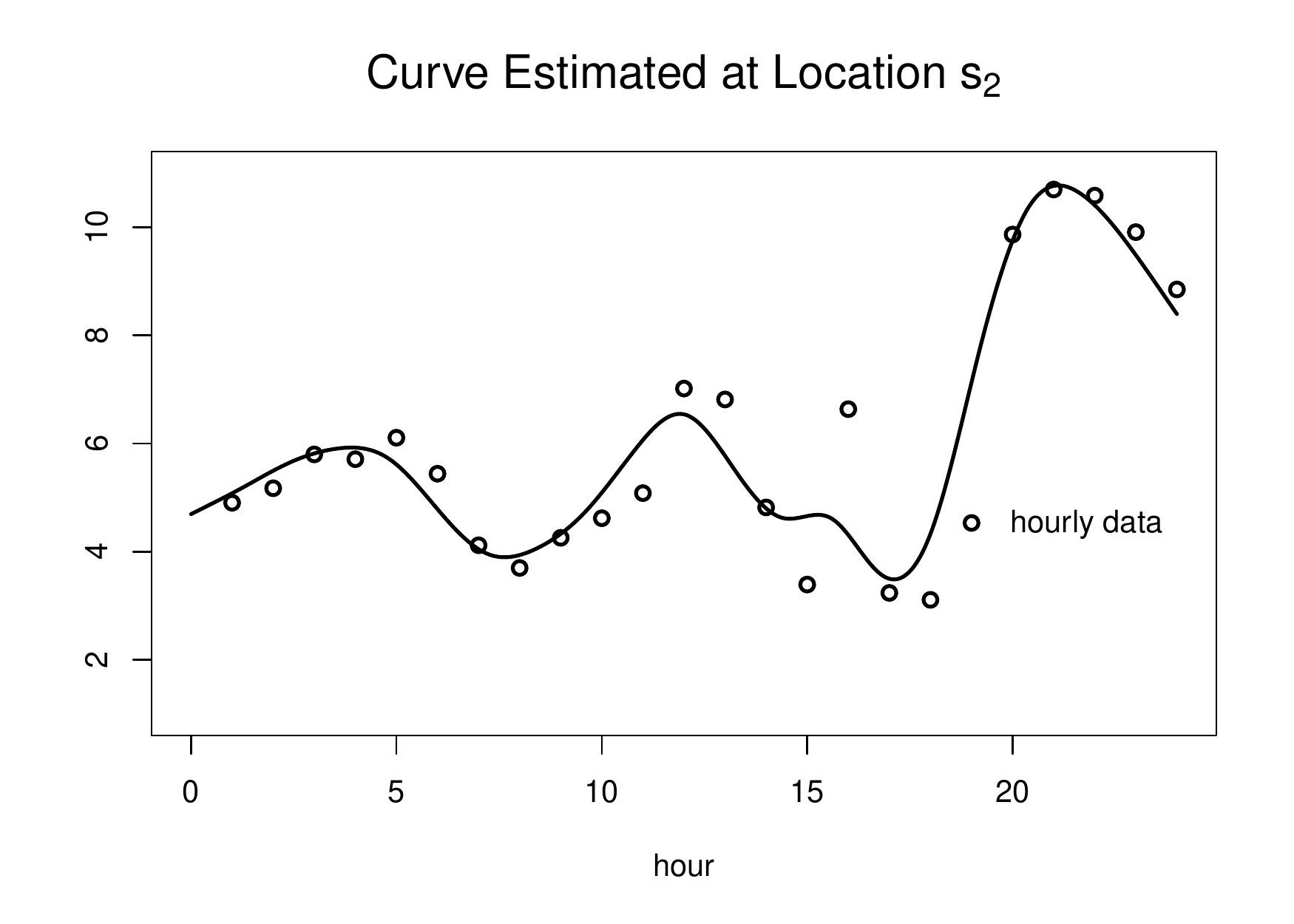} 
\caption{Example of two estimations of the continuous curve at two different locations $\mathbf{s}_{1}$ and $\mathbf{s}_{2}$ separated by $10$ km with latitude fixed. We use $23$ cubic B-spline basis functions.}\label{ExCurveEstimated}
\end{center}
\end{figure}

To illustrate the continuous estimation of the curves, we focus on data observed on June $22$, 2010. For that specific day, we fit $23$ cubic B-spline basis functions for each location, i.e., $\hat{y} (\mathbf{s}_{i};v)= \sum_{k=1}^{23}z_{k}(\mathbf{s}_{i})\eta_{k}(v)$, where $\{\eta_{k}\}_{k=1}^{23}$ are B-spline basis functions. The continuous curves are required to be smooth; thus, we consider a penalization term on the second derivative of $\hat{y} (\mathbf{s}_{i};v)$ with respect to $v$. Specifically, the estimator is such that 
$$y(\mathbf{s}_{i};v_{j}) = \hat{y} (\mathbf{s}_{i};v_{j}) + e_{ij},$$
where the coefficients $z_{k}(\mathbf{s}_{i})$ of $\hat{y} (\mathbf{s}_{i};v)$ are obtained as the solution that minimizes
$$
\sum_{j=1}^{24}  \{ y(\mathbf{s}_{i};v_{j}) - \hat{y} (\mathbf{s}_{i};v_{j}) \}^{2} + \lambda \int  \hat{y}''  (\mathbf{s}_{i};v) \mathrm{d}v,
$$
with the smoothing parameter $\lambda$ fixed.  We select the optimal smoothing parameter by generalized cross-validation through  the $n=529$ locations.  Figure \ref{ExCurveEstimated} shows an example of the curve estimation at two locations. These two locations are separated by $10$ km at the same latitude.

\begin{figure}[!t]
\begin{center}
\includegraphics[scale=.47]{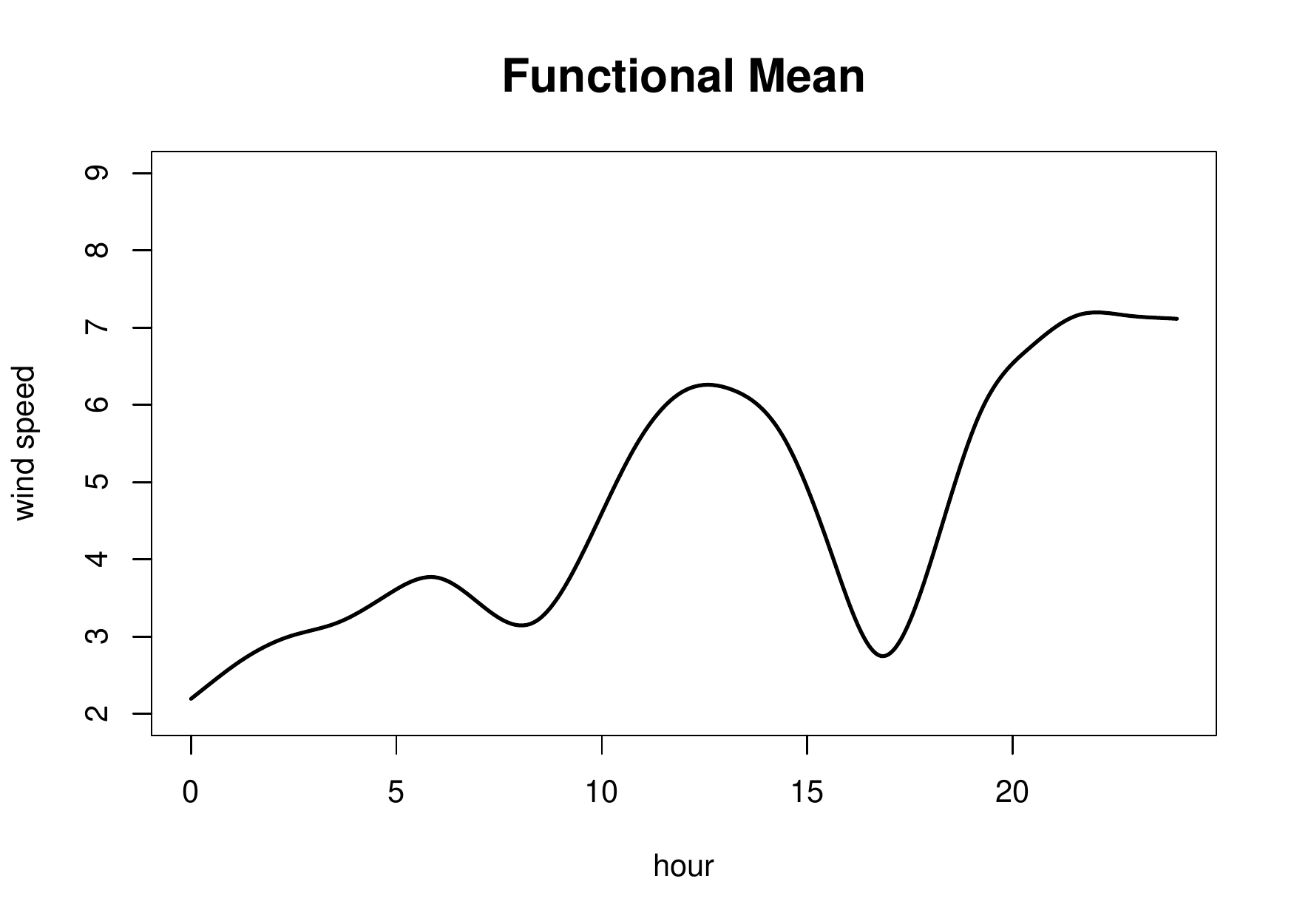} 
\includegraphics[scale=.47]{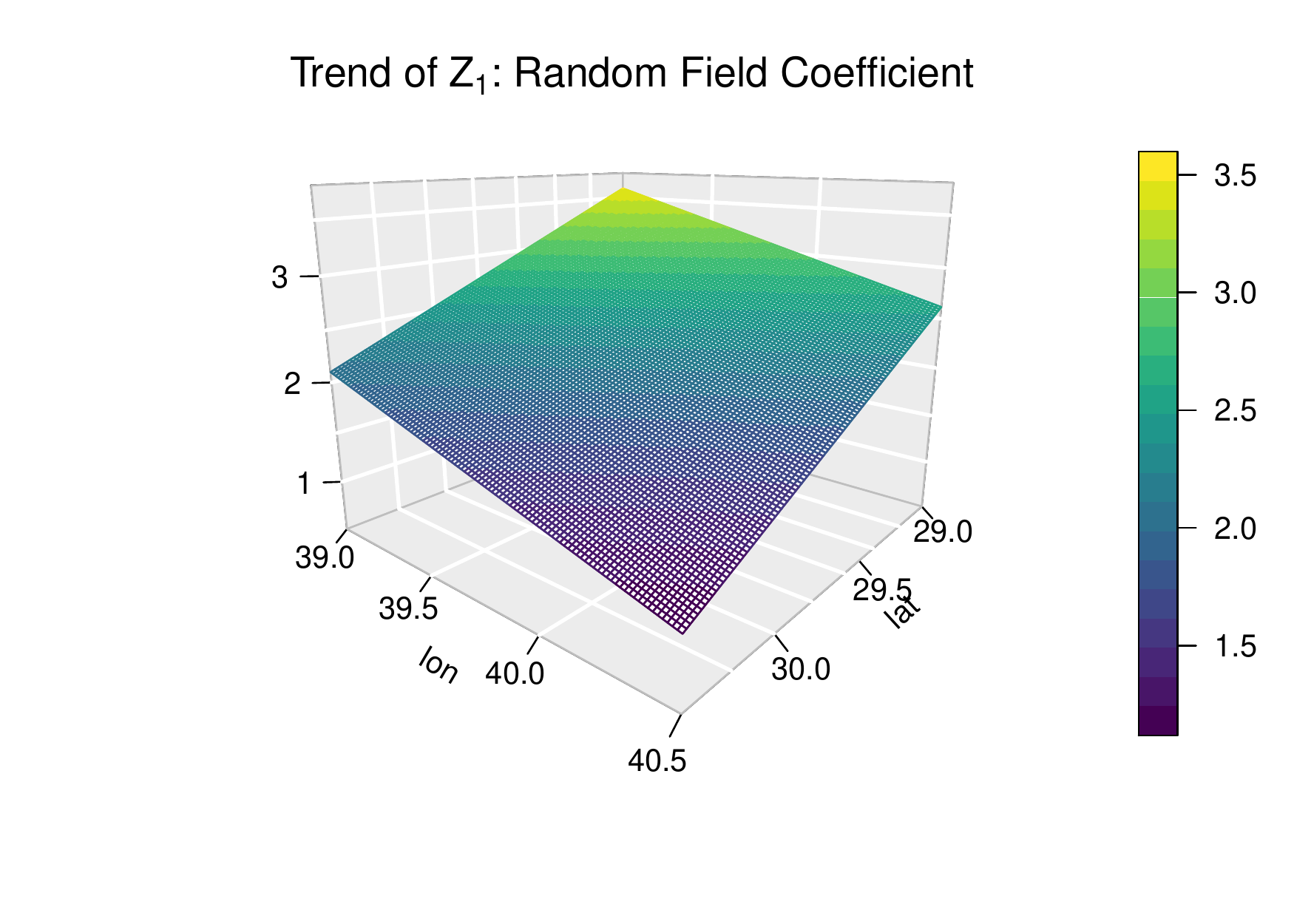} 
\caption{Left: mean curve if the coefficient random fields $\{ Z_{k} ( \mathbf{s}) \}$ are stationary. Right: trend estimation of the first coefficient  random field $Z_{1} ( \mathbf{s} )$ if the coefficient random fields $\{Z_{k} ( \mathbf{s} )\}$ are nonstationary.}\label{WindFdMean}
\end{center}
\end{figure}

Once continuous curves are estimated with basis functions, we can compute the mean (large-scale variation). Since data are measured on a regular grid of points, and if each corresponding coefficient random field $\{ Z_{k} ( \mathbf{s}_{i}) \}_{i=1}^{n}$, with realization $\{  z_{k} ( \mathbf{s}_{i}) \}_{i=1}^{n}$, is stationary, then we can use the empirical mean to estimate each $\mathbb{E}\{ Z_{k} ( \mathbf{s}) \}$. Under this scenario, Figure \ref{WindFdMean} (left) shows the result of the estimated mean curve. 

On the other hand, if the coefficient random fields $\{  Z_{k} ( \mathbf{s}_{i}) \}_{i=1}^{n}$ are nonstationary, then $\mathbb{E}\{ Z_{k} ( \mathbf{s}) \}$ can depend on the location $ \mathbf{s}$. Under this framework, we use the \texttt{surf.gls} command defined in the \textit{spatial} package. Figure \ref{WindFdMean} (right) shows the estimator of $\mathbb{E}\{ Z_{1} ( \mathbf{s}) \}$ for the first random field coefficient $Z_{1} ( \mathbf{s})$. The \texttt{surf.gls} command uses the GLS method to obtain the estimator with an exponential model as the covariance function. From Figure \ref{WindFdMean} (right), we observe some evidence that the mean depends on the location $ \mathbf{s}$.

\begin{figure}[!b]
\begin{center}
\includegraphics[scale=.55]{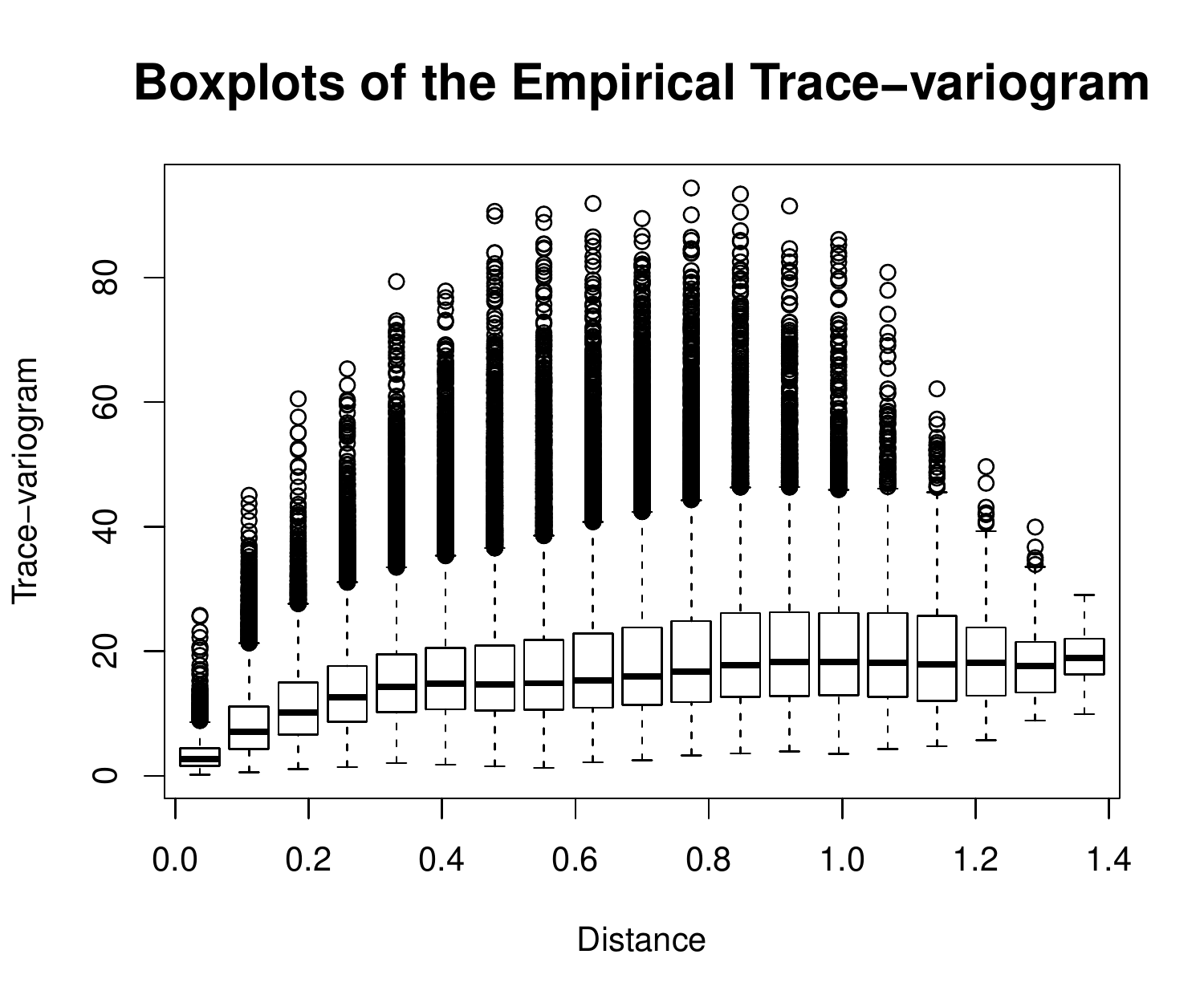} 
\includegraphics[scale=.55]{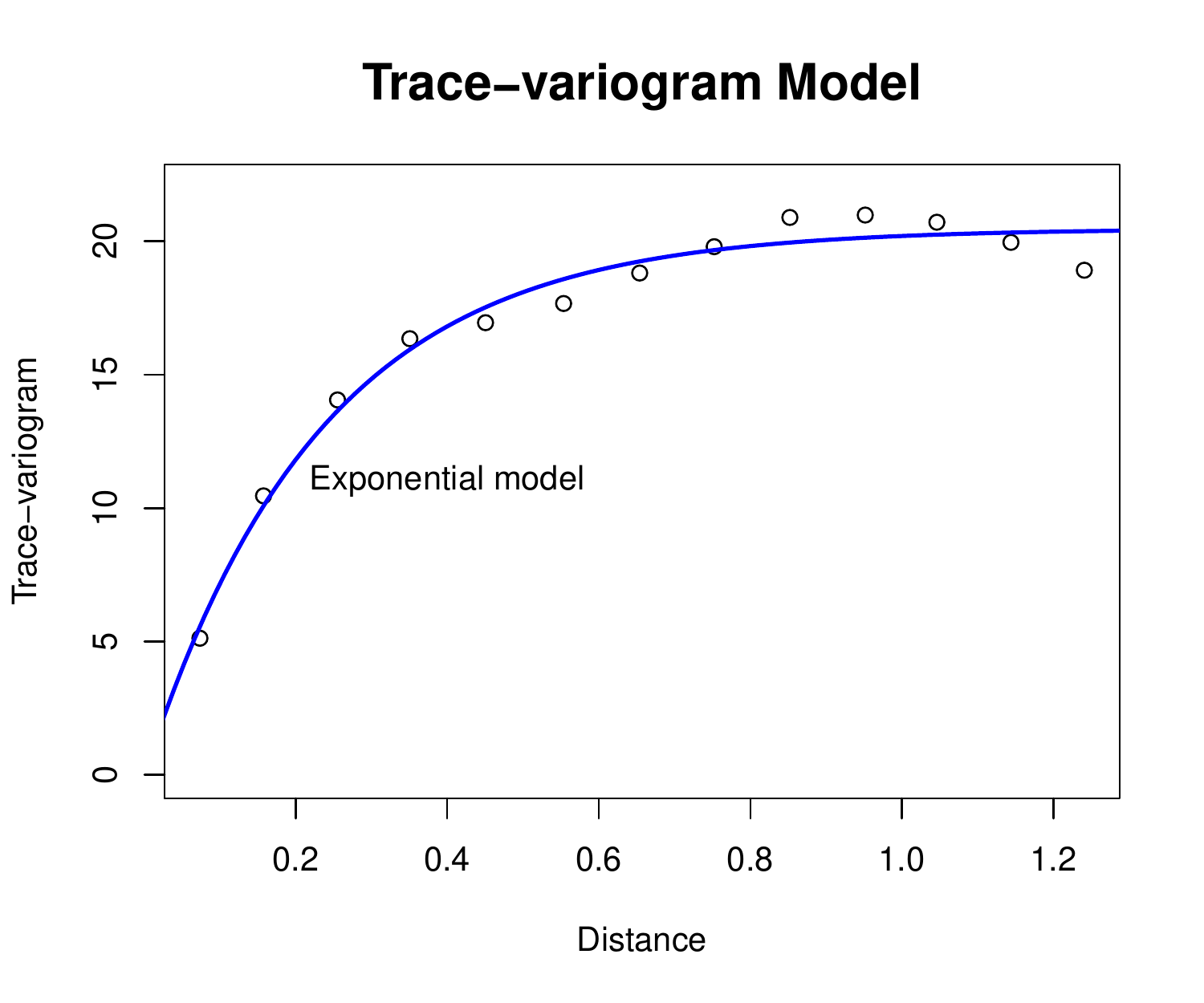}
\caption{Trace-variogram of the daily wind speed. Left: empirical. Right: fitted model (blue).}\label{Wind-Trace-V}
\end{center}
\end{figure}

For the modeling of small-scale variations, we use the trace-variogram defined in \eqref{EmTraVar}.  With the assumption of isotropy, Figure  \ref{Wind-Trace-V} (left) shows boxplots of the corresponding empirical trace-variogram for different distances. The \texttt{fit.tracevariog} command in the \textit{geofd} package allows us to fit four covariance models: spherical, exponential,  Gaussian, and Mat\'ern. In our data, the exponential model is the best model in terms of minimizing the sum of squares errors.  Figure  \ref{Wind-Trace-V} (right) shows the fitted  model. Finally,  ordinary kriging predictors can be obtained with the command \texttt{okfd} in the \textit{geofd} package.

We have used daily wind speed to illustrate general results. Similarly, we can use monthly or yearly wind speed  with $y (\mathbf{s}_{i};v)$ representing the entire year, and repeat the procedure.   

\subsection{Surface Data Approach}

\begin{figure}[!b]
\begin{center}
\includegraphics[scale=.55]{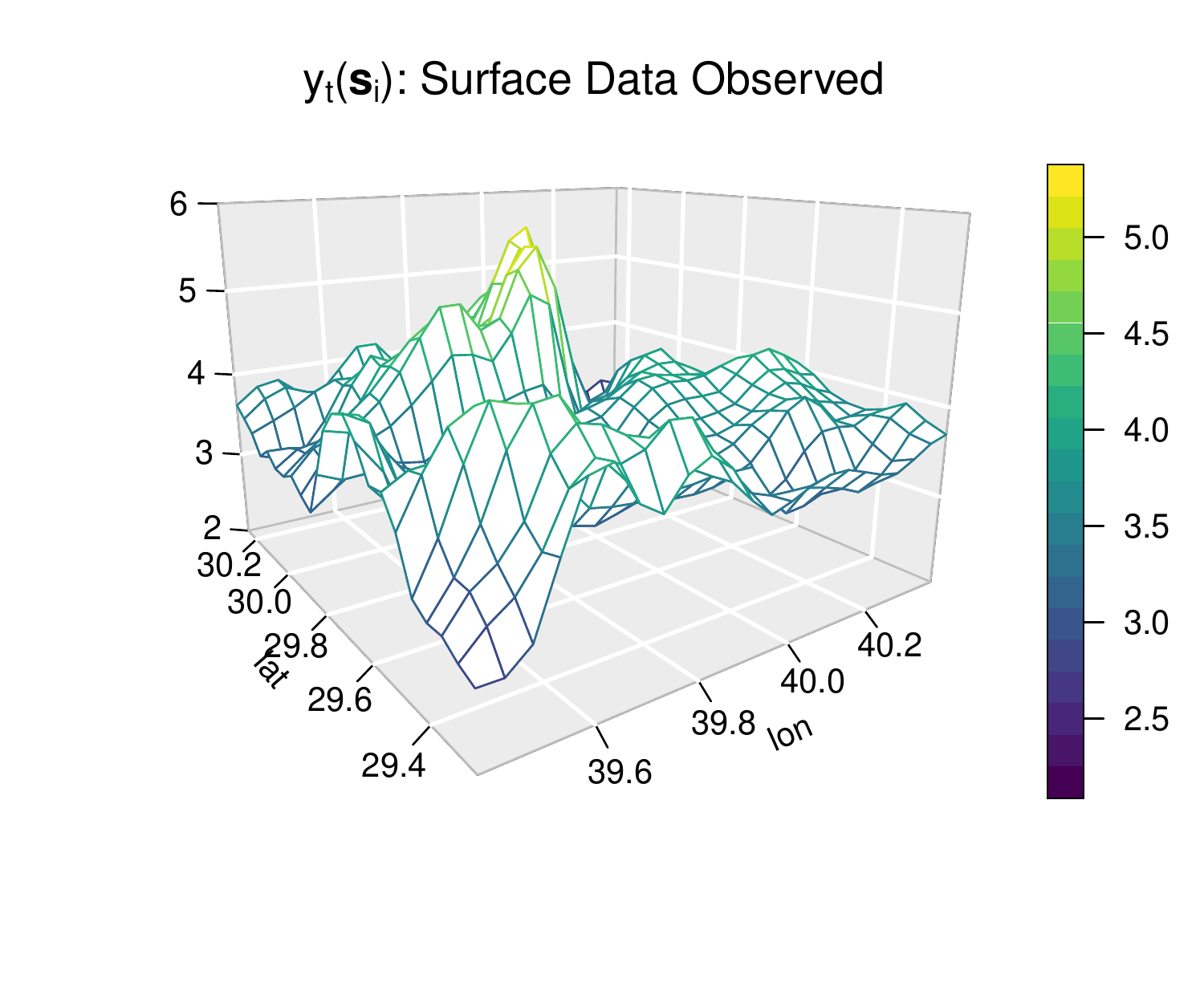} 
\includegraphics[scale=.55]{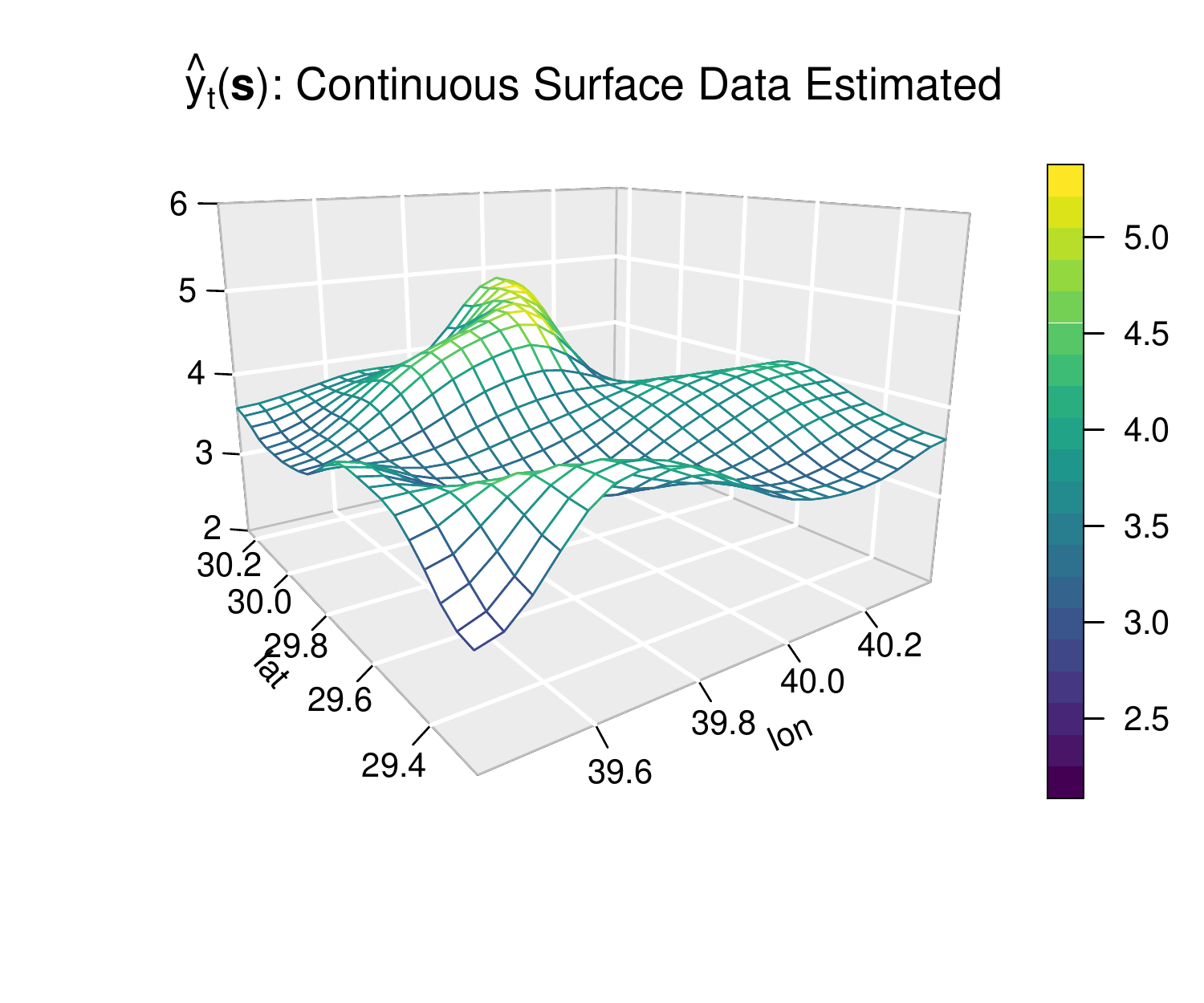} 
\caption{Estimation of surface data for $t=1$. Left: daily average of wind speed on a regular grid of points. Right: the continuous surface estimated with finite elements analysis.}
\label{Wind-surf}
\end{center}
\end{figure}

For illustration purposes, we consider daily average data of the hourly wind speed, that is, $y_{t} (\mathbf{s}_{i})$ represents the daily average wind speed on day $t$, for $t=1,\ldots, 365$ in the year $2010$. 

To estimate the continuous surface $y_{t} (\mathbf{s})$ described in Section   \ref{Est-CS}, we use the \texttt{\small smooth.FEM.basis} 
 command in the \textit{fdaPDE} package. It uses the finite elements analysis and penalizes with the Laplacian. Figure \ref{Wind-surf} shows one surface data for a specific day, $t=1$. In the left panel, we plot the surface data observed on the $115\times 115$ km region centered in Dumah Al Jandal, Saudi Arabia. In the right panel, we plot the estimated continuous surface. The estimated continuous functions $\{\hat{y}_{1} (\mathbf{s}), \ldots, \hat{y}_{365} (\mathbf{s}) \}$ is the surface time series, and can be modeled as in \eqref{SFTS1} and  \eqref{SFTS2}. This approach 
allows us to successfully forecast the next day of surface wind data, $\hat{x}_{T+1} (\mathbf{s})$ and $\hat{y}_{T+1} (\mathbf{s})$, as well as  understanding its temporal dependence. 

Here, we consider a rectangular region domain centered at Dumah Al Jandal, Saudi Arabia, but could also, similarly, 
consider the domain as being the entire kingdom of Saudi Arabia or the  entire world, and study the statistical properties of the surface time series.  In general, the approach of functional data can be used in more complex and large datasets.  

\section{Discussion} \label{Discussion}

In this paper, we have provided an overview of functional data analysis in the case of spatially correlated data. We presented two main approaches,  one when  data are curves observed over space (spatial functional data), and the other when spatio-temporal random fields are considered as a surface (or manifold) time series. These two approaches present a new paradigm of data analysis, in which the continuous processes or random fields are considered as a single entity. Although software packages are still limited for statistical analysis, we believe that this mode of thinking  can be valuable in the context of big data. 

It is a welcomed fact that users can have different tools for data analysis, either with  the surface time series approach or with classical spatio-temporal techniques. The choice will depend on the data, on the phenomena being studied, and on the scientific questions of interest. We appreciate all the effort and work done in both directions. We have attempted to collect all significant references, and we apologize for any unmentioned works.

Similarly, as functional data with temporal dependence arise, we can also consider spatio-temporal functional random fields \citep[][]{BelEtAl2011,LeeEtAl2015} meaning that  a functional time series is observed at each location. For example, $y_{t}(\mathbf{s}_{i};v )$ can represent the functional time series of daily wind speed at location $\mathbf{s}_{i}$, where $t$ is the day index, and $v$ is the time within a day.  Multivariate spatial functional data can also be considered, i.e.,  at each location, we can observe different functional data, such as temperature, precipitation, and humidity.

We conclude that the functional approach opens new areas of research to develop methodologies and theories for the analysis of complex and large spatio-temporal datasets. 

\section*{Acknowledgments} 

We thank Carolina Euan (KAUST), Joydeep Chowdhury (ISI), Laura Sangalli (PoliMi), Soumya Das (KAUST), and Zhuo Qu (KAUST) for comments on this paper. We are grateful to Professor Georgiy Stenchikov’s group, the Atmospheric and Climate Modeling group at KAUST, for producing and providing the high-resolution WRF dataset.


\end{document}